\documentclass[12pt]{dimacs-l}
\usepackage{amssymb}
\usepackage{epsf}

\copyrightinfo{2001}{American Mathematical Society}

\newcounter{FNC}[page]
\def\frankfootnote#1{{\addtocounter{FNC}{2}$^\fnsymbol{FNC}$%
     \let\thefootnote\relax\footnotetext{$^\fnsymbol{FNC}$#1}}}

\newtheorem{thm}{Theorem}[section]
\newtheorem{prop}[thm]{Proposition}

\newtheorem{cor}[thm]{Corollary}
\newtheorem{conj}[thm]{Conjecture}
\newtheorem{ques}[thm]{Question}
\newtheorem{ex}[thm]{Example}
\newenvironment{example}{\begin{ex}\rm}{\end{ex}}

\theoremstyle{remark}
\newtheorem{remark}[thm]{Remark}

\numberwithin{equation}{section}

\newfam\cyrfam
\font\Cyr=wncyr10 at 8pt    
\font\CYR=wncyr10 at 12pt   

\newcommand{\ch}{\mbox{\Cyr ch}}
\newcommand{\CH}{\mbox{\CYR ch}}

\newcommand{\hh}{\hspace{9.3pt}}
\newcommand{\Fdot}{{F\!_\bullet}}
\newcommand{\Edot}{{E_\bullet}}
\newcommand{\New}{{\mbox{\it New}}}
\newcommand{\calA}{{\mathcal A}}
\newcommand{\calM}{{\mathcal M}}
\newcommand{\calS}{{\mathcal S}}
\newcommand{\Gr}{\mbox{\rm Gr}}
\newcommand{\Fln}{{\mathbb{F}\ell_n}}
\newcommand{\Fla}{{\mathbb{F}\ell_{{\bf a}}}}
\newcommand{\sgn}{\mbox{\rm sgn}}
\newcommand{\Char}{\mbox{\rm char}}

\begin{document}

\title[Enumerative Real Algebraic Geometry]{Enumerative Real Algebraic Geometry}

\author{Frank Sottile}
\address{Department of Mathematics\\
        University of Massachusetts\\
        Lederle Graduate Research Tower\\
        Amherst, Massachusetts, 01003\\
        USA}
\email{sottile@math.umass.edu}
\urladdr{http://www.math.umass.edu/\~{}sottile}
\thanks{26 June 2002}
\thanks{Research Supported in part by NSF grant DMS-0070494}
\thanks{Proc. DIMACS workshop on Algorithmic and Quantitative Aspects of
    Real Algebraic Geometry in Mathematics and Computer
       Science, to appear.}

\maketitle

\tableofcontents
\section{Introduction}
Consider the following question.

\begin{ques}\label{q:one}
  Find {\it a priori} information about the number of real solutions to a
  structured system of real polynomial equations
 \begin{equation}\label{eq:polysys}
  0\ =\ f_1\ =\ f_2\ =\ \cdots\ =\ f_N,\qquad\mbox{where each }f_i\in
              {\mathbb R}[x_1,\ldots,x_n]\,. 
 \end{equation}
\end{ques}

Two well-defined classes of structured polynomial systems have been studied
from this point of view---sparse systems, where the structure is encoded by the
monomials in the polynomials $f_i$---and geometric systems, where the
structure comes from geometry.
This second class consists of polynomial formulations of enumerative geometric
problems, and in this case Question~\ref{q:one} is the motivating question of
enumerative real algebraic geometry, the subject of this survey.
Treating both sparse polynomial systems and enumerative geometry
together in the context of Question~\ref{q:one} gives useful insight.

Given a system of polynomial equations~(\ref{eq:polysys}) with $d$ complex
solutions, we 
know the following easy facts about its number $\rho$ of real solutions.
\begin{equation}\label{eq:easy-bounds}
 \begin{minipage}[c]{400pt}
   \begin{enumerate}
    \item[]  We have $\rho\leq d$, \  and   
    \item[]  since $d-\rho$ an even integer, \ 
            ${\displaystyle  \rho\ \geq\ \left\{
              \begin{array}{rl}
                 0&\mbox{ if $d$ is even}\\
                 1&\mbox{ if $d$ is odd\,.}
              \end{array}\right.}$
   \end{enumerate}
 \end{minipage}
\end{equation}
In a surprising number of cases, much better information than this is known.

Structured systems occur in families of systems sharing the same structure.
The common structure determines the number $d$ of complex solutions to a
general member of the family.  
We assume throughout that a general system in any family has
only simple solutions in that each complex solution occurs without
multiplicity. 

Given such a family of structured systems whose general member has 
$d$ complex solutions, perhaps  
the ultimate answer to our motivating question is 
to determine exactly which numbers $\rho$ of real solutions can occur and also
which systems have a given number of real solutions.
Because this level of knowledge may be unattainable, we will be satisfied
with less knowledge.

For example, are the trivial bounds given in~(\ref{eq:easy-bounds}) sharp?
That is, do there exist systems attaining the maximal and minimal number of
real solutions allowed by~(\ref{eq:easy-bounds})?
If these bounds are not sharp, do there exist better sharp bounds?
Perhaps we are unable to determine sharp bounds, but can exhibit systems in a
family with many (or few) real solutions.
This gives lower bounds on the maximum number of real solutions to a system in
a family (or upper bounds on the minimum number).

These answers have two parts: bounds and constructions.
We shall see that bounds (or other limitations) often come from topological 
considerations. 
On the other hand, the constructions often come by deformations from/to a
degenerate situation.
In enumerative geometry, this is the classical technique of special position,
while for sparse systems, it is Viro's method of toric deformations.

In Section 2 we discuss sparse polynomial systems from the point of view of 
Question~\ref{q:one}.
The heart of this survey begins in Section 3, where we discuss some of the
myriad examples of enumerative geometric problems that have been studied from
this perspective. 
In particular, for many enumerative problems the upper
bound~(\ref{eq:easy-bounds}) is sharp.
In Section 4, we concentrate on enumerative problems from the
Schubert calculus, where much work has been done on this question of real
solutions.
Section 5 is devoted to a conjecture of Shapiro and Shapiro, whose study
has led to many recent results in this area.
Finally, in Section 6 we describe new ideas of Eremenko and Gabrielov giving
lower bounds better than~(\ref{eq:easy-bounds}) for some enumerative problems.

\section{Sparse Polynomial Systems}

Perhaps the most obvious structure of a system of $n$ polynomials in $n$
variables
 \begin{equation}\label{eq:sqSystem}
  f_1(x_1,\ldots,x_n)\ =\   f_2(x_1,\ldots,x_n)\ =\ \cdots\ =\ 
  f_n(x_1,\ldots,x_n)\ =\ 0\,,
 \end{equation}
is the list of total degrees of the polynomials $f_1,f_2,\ldots,f_n$.
For such a system, we have the degree or B\'ezout upper bound, which is a
consequence of the refined B\'ezout Theorem of Fulton and 
MacPherson~\cite[\S 1.23]{Fu84a}.

\begin{thm}[B\'ezout Bound] 
 The system~$(\ref{eq:sqSystem})$, where the polynomial $f_i$ has total degree
 $d_i:=\deg(f_i)$, has at most $d_1\cdot d_2\cdots d_n$ isolated complex
 solutions.
\end{thm}

The B\'ezout bound on the number of real solutions is sharp.
For example, if 
 \begin{equation}\label{eq:max-real}
  f_i\ :=\ (x_i-1)(x_i-2)\cdots(x_i-d_i)\quad i=1,\ldots,n\,,
 \end{equation}
then the system~(\ref{eq:sqSystem}) has $d_1\cdot d_2\cdots d_n$ real
solutions.
The reader is invited to construct systems with the minimum possible
number (0 or 1) of real solutions.

A system of polynomial equations with only simple solutions, but with fewer
solutions than the B\'ezout bound is called deficient.
For example, fewer monomials in the polynomials lead to fewer solutions.
We make this idea more precise.
A monomial (or rather its exponent vector) is a point of ${\mathbb N}^n$.
The convex hull of the monomials in $f$ is its {\it Newton polytope},
$\New(f)$, a polytope with integral vertices.
The terms of $f$ are indexed by the lattice points, 
$\New(f)\cap{\mathbb N}^n$, in the Newton polytope of $f$.

\begin{figure}[htb]
$$
  \begin{picture}(141,145)(-24,0)
   \put(  0,  0){\epsfxsize=1.2in\epsfbox{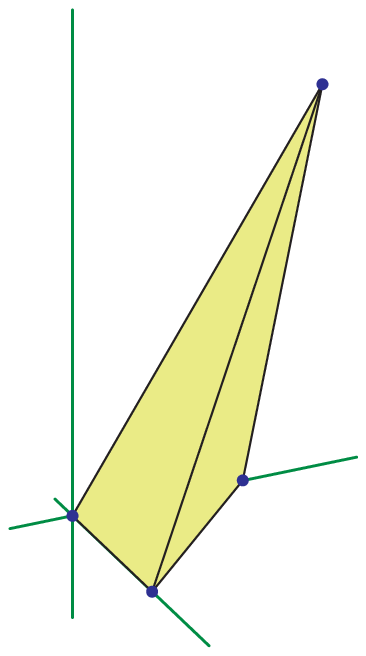}}
   \put(-24, 36){\small $(0,0,0)$}
   \put( 45, 11){\small $(1,0,0)$}
   \put( 60, 32){\small $(0,1,0)$}
   \put( 37,133){\small $(1,1,4)$}
   \put( 76,86){$\New(f)$}
  \end{picture} 
  \hspace{.5in}
  \begin{picture}(125,145)(-24,0)
   \put(  0,  0){\epsfxsize=1.2in\epsfbox{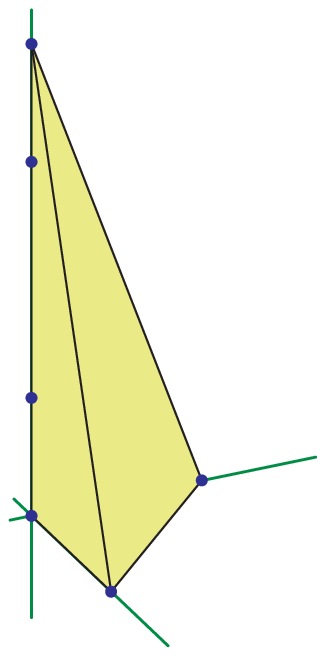}}
   \put(-24, 30){\small $(0,0,0)$}
   \put(-24, 58){\small $(0,0,1)$}
   \put(-24,114){\small $(0,0,3)$}
   \put(-24,142){\small $(0,0,4)$}
   \put( 50,86){$\New(g)$}
   \put( 45, 11){\small $(1,0,0)$}
   \put( 60, 32){\small $(0,1,0)$}
  \end{picture} 
$$
\caption{Newton Polytopes\label{fig:tetrahedrons}}
\end{figure}

Figure~\ref{fig:tetrahedrons} displays the monomials (dots)
and Newton polytopes of the polynomials
 \begin{eqnarray*}
  f &=& 1+ x-y+xyz^4\,,\\
  g &=& 1+ x-y +3z - z^3 + 2z^4\,.
 \end{eqnarray*}

Given polytopes $P_1,P_2,\ldots,P_m\subset{\mathbb R}^n$, their Minkowski sum
is the pointwise sum
$$
  P_1+P_2+\cdots+P_m\ =\ 
       \{ p_1+p_2+\cdots+p_m\mid p_i\in P_i,\ i=1,2,\ldots,m\}\,.
$$
The mixed volume, $\mbox{\it MV}(P_1,P_2,\ldots,P_m)$, of a collection of $m$
polytopes is
$$
  \frac{\partial}{\partial t_1}\frac{\partial}{\partial t_2}
  \cdots\frac{\partial}{\partial t_m}\left[\mbox{Volume of }\;\rule{0pt}{12pt}
  t_1P_1 + t_2P_2 + \cdots + t_mP_m\right]\,.
$$

When there are $n$ equal polytopes,
the mixed volume is $n!\mbox{Vol}(P)$, the normalized volume of the common
Newton polytope $P$, and when each $P_i$ is a line segment of length $d_i$ in
the $i$th coordinate 
(so that $P_i=\New(f_i)$, where $f_i$ is the polynomial given
in~(\ref{eq:max-real})), the mixed volume is the 
B\'ezout number $d_1d_2\cdots d_n$.

Given a list $(P_1,P_2,\ldots,P_n)$ of polytopes in ${\mathbb R}^n$
with vertices in the integral lattice ${\mathbb N}^n$, a {\it sparse
polynomial system}  with this structure is a 
system of polynomials~(\ref{eq:sqSystem}) with $\New(f_i)=P_i$.
These sparse systems may have trivial solutions
where some coordinates vanish.
Thus we only consider solutions in $({\mathbb C}^\times)^n$.
We have the following basic result on the number of solutions to such a sparse
system of polynomials.

\begin{thm}[BKK bound]\label{thm:BKK}
 A sparse polynomial system~$(\ref{eq:sqSystem})$ with structure
 $(P_1,P_2,\ldots,P_n)$ has at most $\mbox{\it MV}(P_1,P_2,\ldots,P_n)$
 isolated solutions in $({\mathbb C}^\times)^n$.
 If the system is generic given its structure, then it has exactly this number
 of solutions in $({\mathbb C}^\times)^n$.
\end{thm}

This result was developed in a series of papers by
Kouchnirenko~\cite{Ko75}, Bernstein~\cite{Bernstein}, and 
Khovanskii~\cite{Kh78}.
For simplicity of exposition, we will largely restrict ourselves to the case
when the polynomials all have the same Newton polytope $P$.

Given a polytope $P$ with
vertices in 
the integral lattice, what are the possible numbers of real solutions to
systems with that Newton polytope?
We shall focus on understanding $\rho(P)$, the maximum number of real
solutions to a sparse system with Newton polytope $P$.
The following example serves as an introduction to this question.

\begin{example}
 Let $d$ be a positive integer, and set
\begin{eqnarray*}
  P_d&:=&\mbox{Convex hull }\{(0,0,0),\ (1,0,0),\, (0,1,0)\,, (1,1,d)\}\,,\ 
  \mbox{ and}\\
  Q_d&:=&\mbox{Convex hull }\{(0,0,0),\ (1,0,0),\, (0,1,0)\,, (0,0,d)\}\,.
\end{eqnarray*}
(Figure~\ref{fig:tetrahedrons} shows $P_4$ and $Q_4$.)
These  tetrahedrons each have normalized volume $d$.

  A general sparse system with Newton polytope $P_d$
\begin{eqnarray*}
 A_1 xyz^d + B_1 x + C_1 y + D_1 &=&0\\
 A_2 xyz^d + B_2 x + C_2 y + D_2 &=&0\\
 A_3 xyz^d + B_3 x + C_3 y + D_3 &=&0
\end{eqnarray*}
is equivalent to a system of the form 
$$
  z^d\ =\ \alpha,\quad x\ =\ \beta,\quad y\ =\ \gamma\,,
$$
where $\alpha,\beta,\gamma\in{\mathbb R}$.  
Thus the original system has $d$ complex solutions, but only 0, 1, or 2 real
solutions, so $\rho(P_d)$ is 1 or 2, depending upon the parity of $d$.

A general sparse system with Newton polytope $Q_d$
is 3 polynomials of the form
$$
  Ax + By + C(z)\,,
$$
where $C(z)$ has degree $d$.
This is equivalent to a system of the form
$$
  f(z)\ =\ 0,\quad x\ =\ g(z),\quad y\ =\ h(z)\,,
$$
where $f(z), g(z)$, and $h(z)$ are real polynomials with
$\deg(f)=d$ but $\deg(g)=\deg(h)=d-1$.
In this case, the general system again has $d$ complex solutions, but there
are systems with any number of real solutions, so $\rho(Q_d)=d$.
\end{example}

\subsection{Polyhedral homotopy algorithm}\label{sec:HS}
The polyhedral homotopy algorithm of Huber and Sturmfels~\cite{HS95}
deforms the sparse system~(\ref{eq:sqSystem})
into a system where the number of solutions is evident.
It gives an effective demonstration of the BKK bound and is 
based upon Sturmfels's generalization~\cite{St94a} of 
Viro's method for constructing real varieties with controlled
topology~\cite{Viro}.

\begin{example}\label{eq:simplex}
Suppose $P$ is a $n$-simplex which meets the integral lattice only at its
vertices.
Translating one vertex to the origin, the others are linearly independent.
(Translating corresponds to division  by a monomial.) 
Let $M$ be the $n\times n$ integral matrix whose columns are these vertices.
Multiplying $M$ by an invertible integral matrix and taking the resulting 
columns to be a basis for ${\mathbb R}^n$ corresponds to a
multiplicative change of coordinates
$$
  (x_1,x_2,\ldots,x_n)\ \longmapsto\ 
  (y^{m_1},y^{m_2},\ldots,y^{m_n})\,,
$$
with $y_1,y_2,\ldots,y_n\in({\mathbb C}^\times)^n$ and $m_1,m_2,\ldots,m_n$
linearily independent integer vectors.
Doing this for the Smith normal form of $M$ 
transforms a polynomial $f$ with Newton polytope $P$ into a polynomial of
the form
 \begin{equation}\label{eq:SNF}
  c_0+c_1 y_1^{d_1}+c_2 y_2^{d_2}+\cdots  +c_n y_n^{d_n}\,,
 \end{equation}
where $d_1\mid d_2\mid \cdots\mid d_n$ and 
$d_1d_2\cdots d_n=n!\mbox{Vol}(P)$, the normalized volume of $P$.
(These $d_i$ are the invariant factors of the integral matrix $M$.)
A system consisting of $n$ general polynomials of the form~(\ref{eq:SNF})
is equivalent to one of the form
$$
  y_1^{d_1}\ =\ \alpha_1\,,\quad y_2^{d_2}\ =\ \alpha_2\,,\quad \ldots\ ,\quad 
  y_2^{d_n}\ =\ \alpha_n\,.
$$
Thus a general system whose Newton polytope is $P$
has $n!\mbox{Vol}(P)$ simple complex solutions.
\end{example}

The combinatorial structure underlying the homotopy algorithm of
Huber and Sturmfels is that of a regular triangulation, which
is a special case of a regular subdivision. 
A regular subdivision $P_w$ of a lattice polytope $P$ is given by a lifting
function 
$$
   w\colon P\cap{\mathbb N}^n\ \longrightarrow\ {\mathbb Q}_{\geq0}
$$
as follows.
Set
$$
  Q\ :=\ \mbox{Convex hull }\{(a,w(a))\mid a\in P\cap{\mathbb N}^n\}\ 
       \subset\ {\mathbb R}^{n+1}
$$
This lifted polytope $Q$ has distinguished lower facets, those facets whose
inward pointing normal vector has positive last coordinate.
Forgetting the last coordinate projects these lower facets into 
${\mathbb R}^n$ (and hence $P$).
Their totality gives the {\it regular subdivision} $P_w$ of $P$.
When all the lattice points in $P\cap{\mathbb N}^n$ lift to vertices of 
lower facets of $Q$, and these lower facets are all simplices, then $P_w$ is a
{\it regular triangulation} of $P$.

\begin{figure}[htb]
$$
   \epsfxsize=1.6in\epsfbox{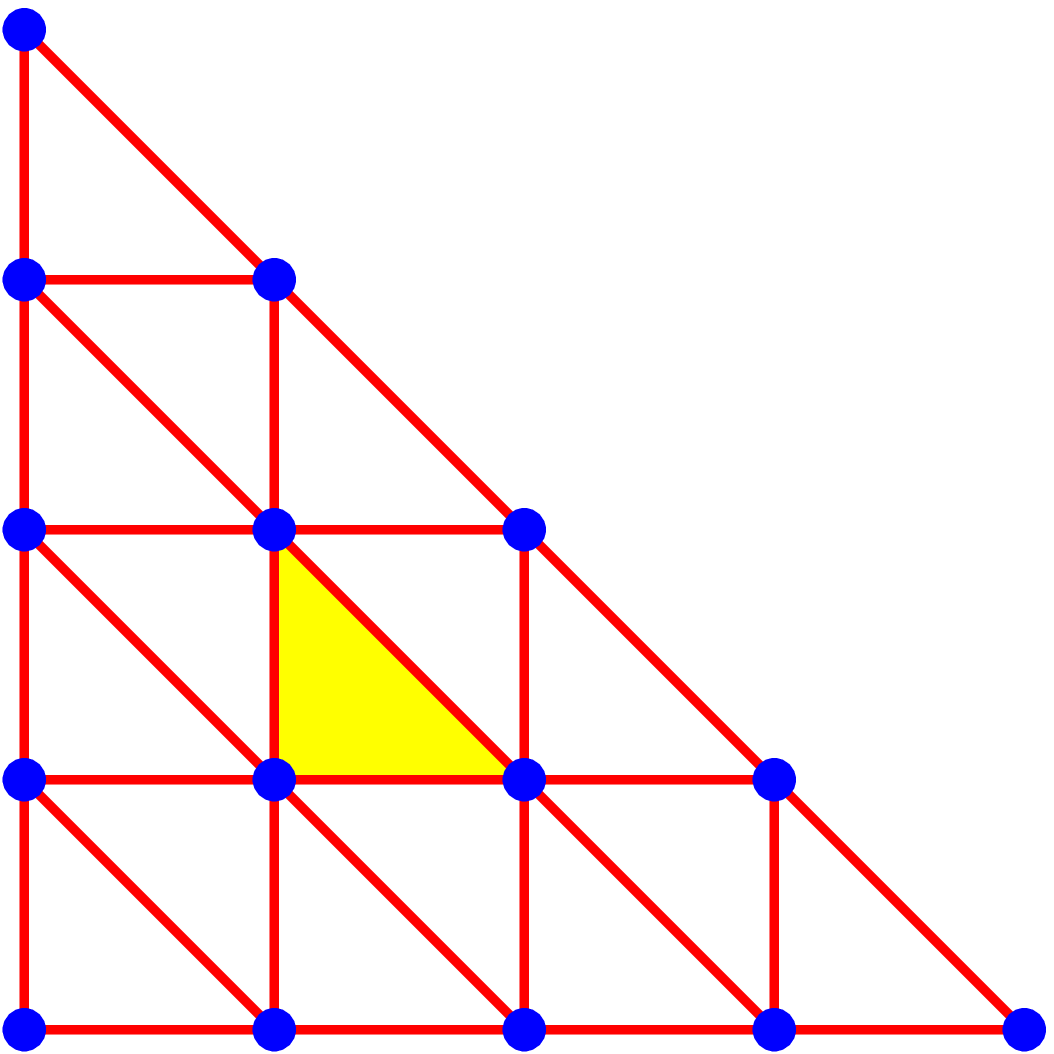}\qquad\qquad
   \epsfxsize=1.6in\epsfbox{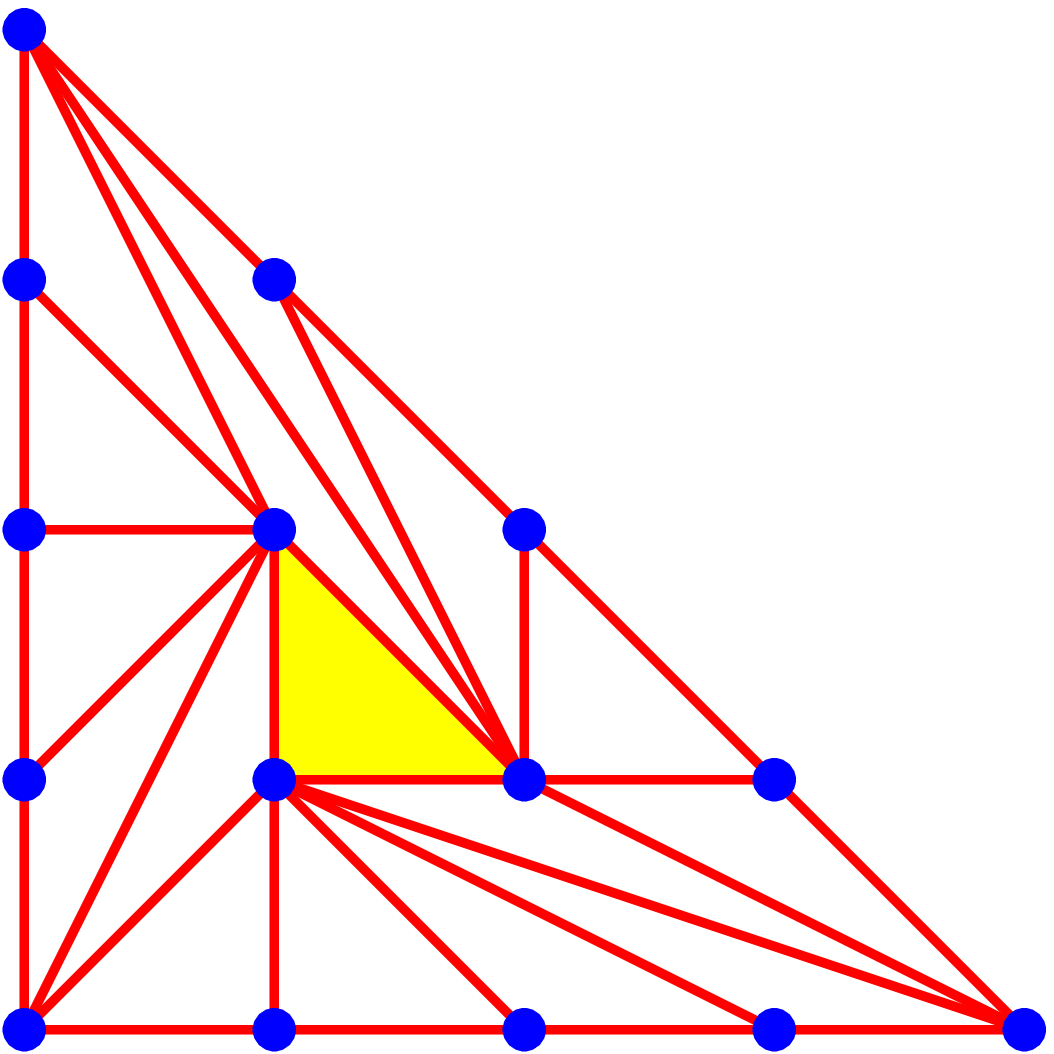}
$$
\caption{Regular and non-regular triangulations\label{fig:regTri}}
\end{figure}

In Figure~\ref{fig:regTri} the triangulation on the left is
regular and 
the triangulation on the right is not regular.
Consider a hypothetical lifting function $w$ for the triangulation on the right.
We assume that $w$ takes the value 0 at the three interior vertices.
The clockwise neighbour of any vertex of the big triangle must be lifted
higher than that vertex. 
(Consider the figure they form with the parallel edge of the
interior triangle.)
Since the edge of the big triangle is lifted to a concave broken path,
this implies that each vertex is lower than its clockwise neighbour,
which is impossible, except in some M.C. Escher woodcuts.

Crucial to this algorithm are polynomial systems derived
from the original system and the  regular subdivision $P_w$ of $P$ by a
lifting function $w$.
Given a polynomial $f$ with Newton polytope $P$ and a face $F$ of the
subdivision $P_w$, consider the sum of terms of $f$ with exponent
vector in $F$, restricting $f$ to the face $F$.
Doing this for each polynomial $f_i$ in our original system, we obtain the
facial subsystem of~(\ref{eq:sqSystem}) given by $F$.

We continue with the algorithm of Huber and Sturmfels.
Given a lifting function $w\colon P\cap{\mathbb N}^n\to{\mathbb Q}_{\geq0}$
and a polynomial $f$ with Newton polytope $P$
$$
    f(x)\ = \sum_{a\in P\cap{\mathbb N}^n} c_a x^a\,,
$$
we multiply the term with exponent $a$ by $t^{w(a)}$ to obtain
$$
   f(x;t)\ :=\ \sum_{a\in P\cap{\mathbb N}^n} c_a x^at^{w(a)}\,.
$$
Modifying all polynomials in the original sparse system in this way gives 
the {\it lifted system}, a family of sparse systems depending upon the
parameter $t$. 
 \begin{equation}\label{eq:spHomotopy}
  f_1(x;t)\ =\ f_2(x;t)\ =\ \cdots\ =\ f_n(x;t)\ =\ 0\,,
 \end{equation}

Solutions to this system are algebraic functions $t\mapsto x(t)$ in the
parameter $t$.
In a neighborhood of 0 in the complex plane, each branch is expressed as a
Puiseaux series
 \begin{eqnarray*}
   x(t)&=&(x_1(t),x_2(t),\ldots,x_n(t))\,,\ \qquad\mbox{with}\\
   x_i(t)&=& t^{u_i}y_i+\ \mbox{higher order terms in $t$}\,,
 \end{eqnarray*}
where the $y_i$ are non-zero constants and the exponents $u_i\in{\mathbb Q}$.
Substituting this expression into~(\ref{eq:spHomotopy}), we obtain
$$
  0\ =\ \sum_{a\in P\cap{\mathbb N}^n} c_a y^a t^{u\cdot a+w(a)}
        +\  \mbox{higher order terms in $t$}\,.
$$
The exponent of $t$ is the dot product $(u,1)\cdot(a,w(a))$.
Thus the terms of lowest order in $t$ correspond to points $(a,w(a))$ in the
lifted polytope where the exponent vector $(u,1)$ achieves its minimum---a
face in the lower hull of the lifted polytope.  
Let $F$ be the corresponding face of the regular subdivision $P_w$.
Removing the common factor $t^{u\cdot a+w(a)}$ from these lowest order
terms gives the restriction of $f$ to the face $F$.
Thus the coefficients $y_i$ of the leading terms of the Puiseaux expansion
are solutions to the facial subsystem of the original
system~(\ref{eq:sqSystem}) given by the face $F$ selected by the initial
Puiseaux exponents.

Suppose the original system is general in the sense that each of its facial
subsystems given by faces $F$ of the regular subdivision has solutions only if
$\mbox{Vol}(F)>0$.
Consider substituting the Puiseaux series into the polynomial
system~(\ref{eq:spHomotopy}) and then taking the limit as $t$ approaches 0.
The resulting polynomial system for the leading coefficients $y$ of the
Puiseaux series are solutions to the corresponding facial system.
By previous assumption, this has no solutions unless $F$ is a facet of the 
subdivision $P_w$.
In that case, the inward pointing normal vector $(u,1)$ to the corresponding
facet of the lifted polytope gives the initial exponent $u$ of the Puiseaux
expansion.

The full Puiseaux series can be reconstructed from its initial terms
(see~\cite{HS95} for details) giving the 1-1 correspondence
$$
  \left\{\begin{minipage}{2.2in}
      Solutions to facial subsystems of~(\ref{eq:sqSystem}) 
      given by facets of $P_w$
  \end{minipage}\right\}\ \Longleftrightarrow\ 
  \left\{\begin{minipage}{2in}
      Branches of the algebraic
       function $t\mapsto x(t)$ near $t=0$
  \end{minipage}\right\}
$$
This is the number of solutions to~(\ref{eq:spHomotopy}) 
for general $t\in{\mathbb C}$, so the number of solutions
to the original sparse system equals the number of solutions to all the facial 
subsystems given by facets of $P_w$.

Now suppose the regular subdivision $P_w$ is a regular triangulation and 
each facial subsystem is general.
Then the facets $F$ of $P_w$ are simplices with no interior vertices.
Since general sparse systems whose Newton polytope is such
a simplex have exactly $n!\mbox{Vol}(F)$ solutions, the number of solutions to
facial subsystems given by facets of $P_w$ is 
exactly the sum of the normalized volumes of all facets of $P_w$, which is
the normalized volume of $P$.

\subsection{Real solutions to sparse polynomial systems}\label{sec:sparse}

How many real solutions are there to a sparse polynomial
system~(\ref{eq:sqSystem}) with common Newton polytope $P$?
By Example~\ref{eq:simplex}, if $P$ is the convex hull of
$v_0,v_1,\ldots,v_n$, a simplex with no interior lattice points, then 
a general system~(\ref{eq:sqSystem}) is equivalent to the system of binomials
 \begin{equation}\label{eq:binomial_system}
   y_1^{d_1}\ =\ \alpha_1,\quad
   y_2^{d_2}\ =\ \alpha_2,\quad\ldots,\quad
   y_n^{d_n}\ =\ \alpha_n\,,
 \end{equation}
with each $a_i\neq 0$.
Here, $d_1\mid d_2\mid \cdots\mid d_n$ are the invariant factors of the
matrix whose $i$th column is $v_i-v_0$ and $d_1d_2\cdots d_n$ is the
normalized volume of $P$.
Following Sturmfels~\cite{St94b}, let $e(P)$ be the number of these
invariant factors which are even.
If $e(P)=0$, so that $d$ is odd, then $P$ is an odd cell.

\begin{prop}\label{prop:simplex}
 The polynomial system~$(\ref{eq:binomial_system})$ has $2^{e(P)}$ real
 solutions if $a_i>0$ whenever $d_i$ is even, and $0$ real solutions otherwise.
 In particular, if $P$ is an odd cell, then there is one real solution.
\end{prop}

\begin{thm}[Sturmfels~{\cite[Corollary 2.3]{St94b}}]
 The maximum number $\rho(P)$ of real solutions to a sparse polynomial
 system~$(\ref{eq:sqSystem})$ with common Newton polytope $P$
 satisfies
$$
   \#\mbox{\rm \ odd cells in }P_w\ \leq\ \rho(P)\,,
$$
 for any regular triangulation $P_w$ of $P$.
\end{thm}

\begin{proof}
 In the limit as $t\to 0$, the lifted system~(\ref{eq:spHomotopy}) given by a
 regular triangulation $P_w$ of $P$ becomes the disjunction of facial subsystems
 of~(\ref{eq:sqSystem}), one for each facet $F$ of $P_w$.
 Thus the number
 of real solutions in the limit is a constructive lower bound for $\rho(P)$.

 By Proposition~\ref{prop:simplex}, the number of odd cells in $P_w$ is a
 lower bound on the number of real solutions to the facial subsystem given by
 the facets of $P_w$.
\end{proof}

Proposition~\ref{prop:simplex} also gives an upper bound on the
limiting number of real solutions $\rho$ of the lifted
system~(\ref{eq:spHomotopy}) as $t\to 0$ (due to 
Sturmfels~\cite[Theorem 2.2]{St94b}.)
 \begin{equation}\label{eq:limit-bound}
  \rho\ \leq\ \sum_{F\ \mbox{\scriptsize a facet of }P_w} 2^{e(F)}\,.
 \end{equation}
More sophisticated accounting of the possible signs of the coefficients of
facial subsystems improves this bound.
This accounting is accomplished in~\cite{PS96} and~\cite{IR96}, leading to a
combinatorial upper bound for such limiting systems.
Itenberg and Roy~\cite{IR96} show there is a system~(\ref{eq:sqSystem}) for
which this upper bound is attained, and thus obtain
$$
  \mbox{Combinatorial upper bound for limiting systems \ }\leq\ \rho(P)\,.
$$
Itenberg and Roy then conjectured that this combinatorial bound was in
fact the global bound, that is, they conjectured that the maximal number of
real solutions occurs in limiting systems.
They also gave a similar bound for $\rho_+(P)$, the number of solutions with
positive coordinates.
This was too optimistic, for Li and Wang~\cite{LW98} found a
remarkably simple counterexample to this conjecture of Itenberg and
Roy\frankfootnote{Strictly speaking, this counterexample is to their more general
conjecture concerning {\it mixed systems} which are systems of polynomials with
possibly different Newton polytopes.}
$$
        y-x-1\ =\ 200-100y^3+900x^3-x^3y^3\ =\ 0\,.
$$
This system has 3 positive solutions
$$
     (0.317659, 1.317659),\quad 
     (.659995, 1.659995),\quad\mbox{and}\quad
     (8.12058, 9.12058)\,,
$$
but the combinatorial upper bound is 2.
Thus systems with the maximal number of real solutions cannot  in general be
constructed with these limiting techniques.

We still have more questions than answers.
Among the questions are:
\begin{enumerate}
 \item
       Improve this lower bound for $\rho(P)$ (or for $\rho_+(P)$) of Itenberg
       and Roy. 
  \item
        Find new general methods to construct systems with many real solutions.
  \item 
        Determine which polytopes $P$ have the property that
        $\rho(P)<n!\mbox{Vol}(P)$, that is, not all the solutions can be real.
\end{enumerate}

\subsection{Fewnomial bounds}\label{sec:fewnomial}
These definitions, constructions, and results apply to
what have come to be known as fewnomial systems.
A fewnomial is a polynomial $f$ with {\it few} mo{\it nomials}---the
monomials of $f$ are members of some set $\calA$ not
necessarily equal to the lattice points within its convex hull.
In particular, the results of Section~\ref{sec:sparse} give lower bounds on
the maximum number of real solutions $\rho(\calA)$ to a system of fewnomials
whose monomials are from $\calA$.
Here, we use a regular triangulation $\Delta_w$ of the point set $\calA$
induced from a lifting function $w\colon\calA\to{\mathbb Q}$.

When $n=1$, consider the binomial (a fewnomial) system
$$ 
  x^d-1\ =\ 0\,.
$$
This has $d$ complex solutions.
Similarly, the number of complex solutions to a general fewnomial system
is equal to its BKK bound.
The above binomial system has either 1 or 2 real solutions, and so we see that
the number of real solutions to a fewnomial system should be less than its BKK 
bound. 
The question is: How much less?

Khovanskii~\cite{Kh78,Kh91} established a very general result concerning
systems where each $f_i$ is a {\it polynomial} function of the monomials $x^a$
for $a\in\calA$.
He proves that the number of real solutions to such a system is at most
$  2^n2^{\binom{N}{2}}(n+1)^{N}$, where $N=\#\calA$, the number of monomials
in $\calA$.
When the polynomial functions are linear, they are polynomials with monomials
from $\calA$, and hence we have Khovanskii's fewnomial bound.
$$ 
  \rho(\calA)\ \leq \ 2^n\cdot 2^{\binom{N}{2}}(n+1)^N\,.
$$
While this bound seems outrageously large, it does not depend upon the
volume of the convex hull of $\calA$, but rather on the algebraic
complexity---the ambient dimension $n$ and the size $N$ of $\calA$.
That such a bound exists at all was revolutionary.

We compare this complexity bound to the combinatorial upper
bound~(\ref{eq:limit-bound}) on the 
number of real solutions to a lifted fewnomial system~(\ref{eq:spHomotopy}) in
the limit as $t\to 0$.
The invariant $e(F)$ of a facet of the regular triangulation $\Delta_w$ of
$\calA$ is at most $n$.
Thus
$$
  \rho\ \leq\ 2^n\cdot \#\mbox{facets of }\Delta_w\ \left(\leq\ 
              2^n\cdot \binom{N}{n+1}\right)\,,
$$
as a facet of $\Delta_w$ involves $n+1$ points
of $\calA$.
This bound is typically much lower than Khovanskii's bound.
For example, consider two  trinomials in two variables.
Here $n=2$, and after multiplying each equation by a suitable monomial, we
have $\#\calA = 5$.
Thus we have Khovanskii's fewnomial bound
$$
  \rho(\calA)\ \leq\ 2^2\cdot 2^{\binom{5}{2}}\cdot 3^5\ =\ 995,328\,.
$$
In contrast, a triangulation of 5 points in the plane has at most 5 simplices
$$
   \epsfxsize=1.2in\epsfbox{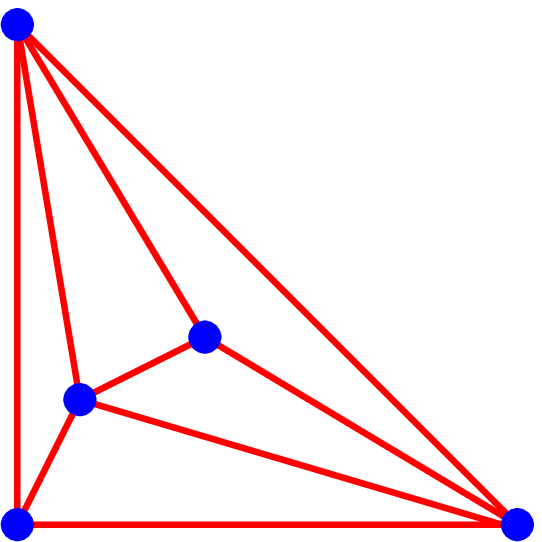}
$$
and so the bound $\rho$ for limiting lifted systems is
$$
  \rho\ \leq\ 2^2 \cdot 5 \ =\ 20\,.
$$

It may be more feasible to look for bounds on
$\rho_+(\calA)$, the number of solutions with positive coordinates.
Kouchnirenko made\frankfootnote{Apparently, Kouchnirenko did not believe this
bound. Nonetheless, the conjecture and its attribution have 
passed into folklore.}
the following conjecture.

\begin{conj}[Kouchnirenko {(see~\cite[p.~300]{BCS97})}]\label{conj:Kouch}
 A system of $n$ polynomials in $n$ variables with only
 simple solutions whose $i$th polynomial has $m_i$ monomials has at most
$$
  (m_1-1)(m_2-1)\cdots(m_n-1)
$$
solutions with positive coordinates.
\end{conj}

For a system of two trinomials in 2 variables, Kouchnirenko's conjecture asserts
that  $\rho_+\leq 4$.
In 2000, Haas~\cite{Ha00} gave the counterexample
 \begin{eqnarray*}
  x^{108}+1.1y^{54}-1.1y&=&0\\
  y^{108}+1.1x^{54}-1.1x&=&0,
 \end{eqnarray*}
a system of two trinomials with 5 solutions having positive
coordinates. 
Although Kouchnirenko's conjecture is false, the question remains:
Is the true value for $\rho_+(\calA)$ (or $\rho(\calA)$) closer to
Khovanskii's bound or to the number in Kouchnirenko's conjecture?
Recent work suggests that it is the latter.

\begin{thm}[Li, Rojas, and Wang~\cite{LRW01}]\label{thm:lrw}
 A system of $2$ polynomials in $2$ variables
$$
   f_1(x,y)\ =\ f_2(x,y)\ =\ 0\,,
$$
 where $f_1$ has $3$ terms and $f_2$ has $m$ terms, and every solution is
 simple, has at most
$$
   2\cdot (2^{m-1}-1)\,,
$$
 solutions with positive coordinates.
\end{thm}

When $m=3$, so we have two trinomials in 2 variables, this gives a bound of 6.
There is presently no known example with more that 5 positive solutions.
Does Theorem~\ref{thm:lrw} give the correct value of 6, or is Haas's
construction with 5 solutions the best possible?

\begin{remark}
 Khovanskii's bound for the number of real solutions with positive
 coordinates holds also for systems of power
 functions, 
 where the exponents of monomials are arbitrary real numbers.
 One might suppose that this added generality is the source of the large size
 of his bound and its apparent lack of tightness for polynomial systems.
 However, Napoletani~\cite{Na01} has shown that the complexity bounds are the
 same for both polynomials and for power functions.
\end{remark}

\section{Enumerative Real Algebraic Geometry}
 
In his treatise on enumerative geometry, Schubert~\cite{Sch1879} declared
enumerative geometry to be concerned with all questions of the following form:
How many geometric figures of a fixed type satisfy certain conditions?
This includes problems as diverse as the number of lines on a cubic surface
(27) and the number of lines meeting four fixed lines in 3-space (see
Section~\ref{sec:SchubertCalculus} for the answer).
These are archetypes for 
two distinct classes of enumerative geometric problems.
For the purpose of this survey, we ignore the first class (except in 
Section~\ref{old-real}) and concentrate on the second.
Specifically, we consider conditions imposed by geometric figures that
may be moved independently into general position.

Enumerative geometry (in the broad sense) had a great
flourishing in the 19th century in the hands of Chasles, Schubert, Halphen,
Zeuthen, and others.
(The survey of Kleiman~\cite{Kl76} is a good historical source.)
At that time, it had long been recognized that it was necessary to work over
the complex numbers to ensure the existence of solutions.
We know of only a handful of cases where the number of real solutions was
considered.
(We discuss some in Section~\ref{old-real}.)
Asking how many solutions can be real is the motivating question of
enumerative real algebraic geometry and an analog of Question~\ref{q:one}.
 
\begin{ques}\label{ques:reality}
 In a given enumerative geometric problem, if the general figures are chosen
 to be real, how many of the solution figures can be real?
\end{ques}

For example, how many of the lines on a real cubic surface can be real?
(Answer: all 27.) How many of the lines meeting four given real lines can be
real? (Answer: all can.)
This question was posed by Fulton~\cite[p.~55]{Fu96b}:
``The question of how many solutions of real equations can be real is still
very much open, particularly for enumerative problems.''

This problem is fundamentally hard.  
Of the geometric figures satisfying real conditions, some will be real while
the rest will occur in complex conjugate pairs, and the number which are
real will depend subtly on the configuration of the conditions.
Despite this difficulty, this is an important question often asked in
applications.

One fruitful variant has been whether it is possible that {\it all} solution
figures can be real.
We call an enumerative problem {\it fully real} if there are general real
conditions for which all solution figures are real.
That is, if the upper bound of~(\ref{eq:easy-bounds}) is sharp.
In light of the situation for sparse polynomial systems, it is quite
surprising that there are no known enumerative problems which are not fully
real. 
For this it is important that the conditions are imposed by general figures.
A related question is whether the opposite situation can occur:
Is it possible to have no (or only one) real solutions?
We give some examples in Section~\ref{sec:lagriangian}.
We shall see in Section~\ref{sec:lower} that there are many enumerative
problems whose number of real solutions is always at least 2.

In the above passage, Fulton~\cite{Fu96b} goes on to ask:
``For example, how many of the 3264 conics tangent to five general conics can
be real?''
He answered this question in the affirmative in 1986, but did not publish that
result.
Later, Ronga, Tognoli, and Vust~\cite{RTV97} gave a careful argument that all
3264 can be real.
This example is very striking, both for the number, 3264, and because this
problem of conics has long been an important testing ground for ideas in
enumerative geometry.

One difficulty with enumerative real algebraic geometry is that there are few
techniques or theorems with a wide range of applicability.
Sometimes a direct calculation suffices (Sections~\ref{sec:ratcubics}
and~\ref{sec:tangent}) or more commonly real solutions are constructed by
deforming from a special configuration, as in the homotopy algorithm of Huber
and Sturmfels in Section~\ref{sec:HS}.
Attempts to formalize this method include  `Schubert
induction'~\cite{So00c}; this is presented in the proof of
Theorem~\ref{thm:real-trans}, and used to establish the other results of
Section~\ref{sec:SchubertCalculus}.
Another formalization is the notion of `real effective algebraic
equivalence'~\cite{So97b}. 
That (together with Theorem~\ref{thm:special}) gives such results
as~\cite[Theorem~18]{So97b}: The enumerative problem of
$\binom{2n-2}{n-1}n^{2n-3}$ codimension 2 planes meeting $2n-2$ rational
normal curves is fully real and~\cite[Corollary 5]{So97c}: The enumerative
problem of 11,010,048 2-planes in ${\mathbb P}^5$ meeting 9 general Veronese
surfaces is fully real.

We devote the rest of this section to a description of some enumerative
problems in the context of Question~\ref{ques:reality}.

\subsection{Enumerative problems not involving general conditions}
\label{old-real}
We give a brief tour of some enumerative problems that do not involve general
conditions, but nonetheless raise some interesting questions regarding real
solutions. 

Schl\"afli (see the survey of Coxeter~\cite{Co83}) showed there are 4
possibilities for the number of real 
lines on a real cubic surface: 27, 15, 7, or 3.
Recent developments in enumerative geometry (mirror symmetry) have led to the
solution of a large class of similar enumerative problems involving, among
other things, the number of rational curves of a fixed degree on a Calabi-Yau
threefold.
(See the book of Cox and Katz~\cite{CK00}.)
For example, on a general quintic hypersurface in ${\mathbb P}^4$ there are 
2875 lines~\cite{Ha79}, 609,250 conics~\cite{Ka86}, and 371,206,375 twisted
cubics. 
The number of twisted cubics and higher degree rational curves was computed in
the seminal paper of Candelas, de la Ossa, Green, and Parkes~\cite{COGP92}.
How many of the curves can be real in problems of this type?
For example, how many real lines can there be on a real quintic hypersurface
in ${\mathbb P}^4$?\smallskip

A real homogeneous polynomial $f(x)$ is positive semi-definite (psd) if
$f(x)\geq 0$ whenever $x$ is real.
Hilbert~\cite{Hi1888} proved that a psd ternary quartic is a sum of three
squares of real quadratic forms.
In fact, a general quartic is a sum of three squares of {\it complex} quadratic
forms in 63 different ways~\cite{Wa91}.
Powers and Reznick~\cite{PR00} studied the question of how many ways one may
represent a ternary quartic as a real sum or difference of three squares.
In every instance, they found that 15 of the 63 ways involved
real quadratic forms.
Is it true that a general psd quartic is a sum or a differrence of three 
real squares in exactly 15 different ways?
In how many ways is it a sum of three squares?\frankfootnote{With Powers, 
Reznick, and Scheiderer, we have shown that the answers to these two questions are
15 and 8, respectively~\cite{PRSS}.}\smallskip

\begin{figure}[htb]
$$
   \epsfxsize=1.35in\epsfbox{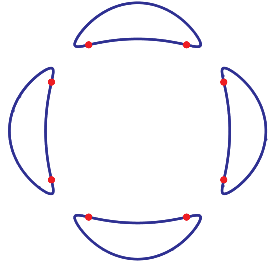}
$$
\caption{Hilbert's quartic: a plane quartic with 8 flexes\label{fig:Hilbert}}
\end{figure}

A general plane curve $C$ of degree $d$ has $3d(d-2)$ flexes.
These are the points on $C$ where the Hessian determinant of the form defining
$C$ vanishes.
Since the Hessian determinant has degree
$3(d-2)$, we expect there to be $3d(d-2)$ such points.
This involves intersecting the curve with its Hessian curve, and 
{\it not} with a general curve of degree $3(d-2)$.
A real smooth plane cubic has 3 of its 9 flexes real.
Zeuthen~\cite{Ze1874} found that at most 8 of
the 24 flexes of a real plane quartic can be real.
An example of a plane quartic with 8 real flexes is provided by the Hilbert
quartic~\cite{Hi1891}, which is defined by
$$
  (x^2+2y^2-z^2)(2x^2+y^2-z^2) + z^4/100\ =\ 0\,.
$$
We display this curve in Figure~\ref{fig:Hilbert}, marking the flexes with
dots. 
Klein~\cite{Klein} later showed that a general real plane curve 
has at most 1/3 of its flexes real.\smallskip

Harnack~\cite{Ha1876} proved that a smooth real algebraic curve of genus $g$
has at most $g+1$ topological components, 
and he constructed real algebraic curves of genus $g$ with $g+1$ components.
In particular, a plane curve of degree $d$ has genus $g=(d-1)(d-2)/2$ and
there are real plane curves of degree $d$ with $g+1$ components.
(An example is provided by Hilbert's quartic, which has genus
$3$.)
Finer topological questions than enumerating the components leads to (part of)
Hilbert's 16th problem~\cite{Hi1902}, which asks for the determination of the
topological types of smooth projectively embedded real algebraic
varieties.\smallskip 

\begin{figure}[htb]
$$
   \epsfxsize=1.1in\epsfbox{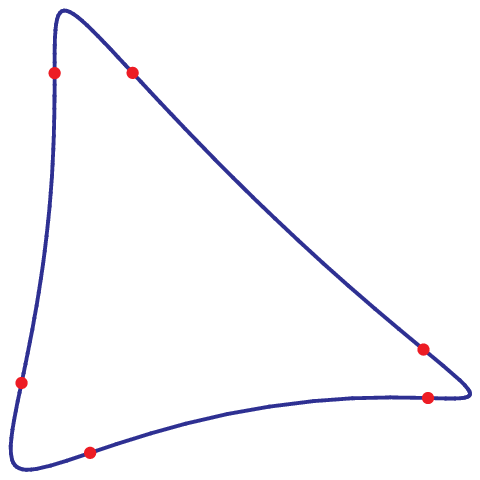}\qquad\qquad
   \epsfxsize=1.1in\epsfbox{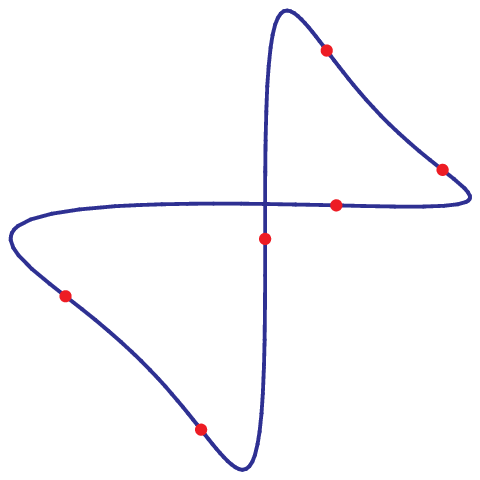}
$$
\caption{Rational quartics with 6 real flexes\label{fig:inflected}}
\end{figure}

A variant concerns rational plane curves of degree $d$.
A general rational plane curve of degree $d$ has $3(d-2)$ flexes and
$g=(d-1)(d-2)/2$ nodes.
Theorem~\ref{thm:reality} implies that there exist real rational plane curves
of degree $d$ with all $3(d-2)$ flexes real, which we call maximally inflected
curves.
See Section~\ref{eq:ratreal} for the connection.
Such curves have at most $g-d+2$ of their nodes real, and there exist
curves with the extreme values of $0$ and of $g-d+2$ real nodes~\cite{KhS}.
For example, a rational quartic ($d=4$) has 6 flexes and $g=3$ nodes.
If all 6 flexes are real, then at most one node is real.
Figure~\ref{fig:inflected} shows maximally inflected quartics with 
and 0 and 1 nodes.
The flexes are marked by dots.\smallskip

Recently, Huisman asked and answered a new question about real curves.
A component $X$ of a real algebraic curve is a {\it psuedoline} if its homology
class $[X]$ in  $H^1({\mathbb P}^n_{\mathbb R},{\mathbb Z}/2{\mathbb Z})$ is
non-zero, and an {\it oval} otherwise.

\begin{ques}
 Given a smooth (irreducible over ${\mathbb C}$) real algebraic curve $C$ in
 ${\mathbb P}^n$ of genus $g$ and degree $c$, how many real hypersurfaces of
 degree $d$ are tangent to at least $s$ components of $C$ with order of
 tangency at least $m$?
\end{ques}

Let $\nu$ be that number, when there are finitely many such hypersurfaces.

\begin{thm}[Huisman~{\cite[Theorem 3.1]{Hu01}}]
 When $s=g$ and $gm=cd$, and the restriction
$$
  H^0({\mathbb P}^n,{\mathcal O}(d))\ \longrightarrow\ 
  H^0(C,{\mathcal O}(d))
$$
 is an isomorphism, then $\nu$ is finite.
 Moreover, $\nu$ is non-zero if and only if
\begin{enumerate}
 \item Both $m$ and $d$ are odd and $C$ consists of exactly $g$ psuedolines,
   or
 \item $m$ is even and either $d$ is even or all components of $C$ are ovals.
\end{enumerate}
In case {\rm (1)}, $\nu=m^g$, and in case {\rm (2)},
$$
  \nu\ =\ \left\{\begin{array}{lcl}
             (g+1)m^g&\ &\mbox{if $C$ has $g+1$ components}\\
                  m^g&\ &\mbox{if $C$ has $g$ components.}
                \end{array}\right.
$$
\end{thm}

It is notable that this problem can only be stated over the real numbers.

\subsection{The Stewart-Gough platform}
\begin{figure}[htb]
$$
  \setlength{\unitlength}{1.1pt}
  \setlength{\unitlength}{1.3pt}
  \begin{picture}(269,200)(-7,0)
   \put(  9, 0){\epsfxsize=4.018in\epsfbox{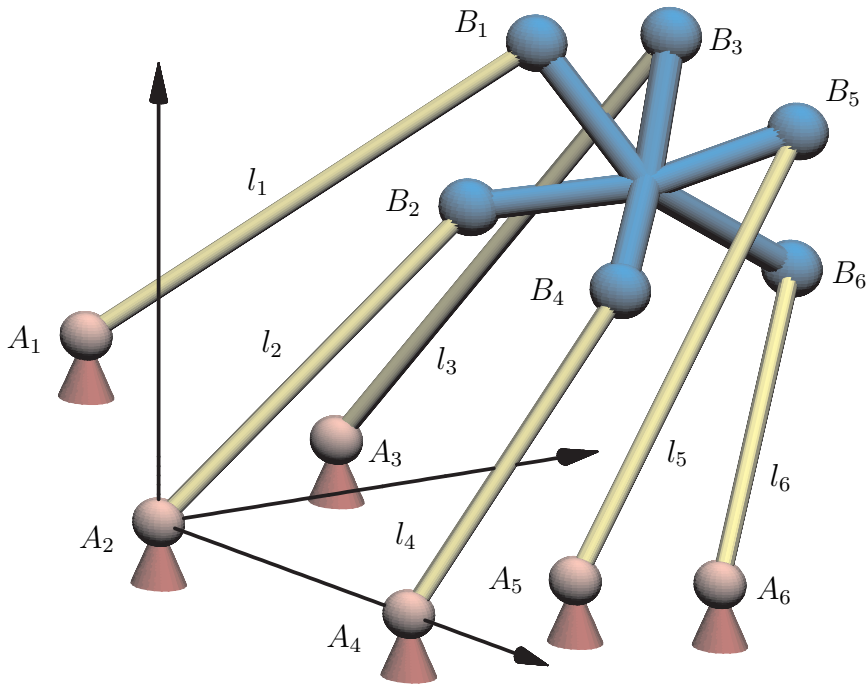}}
   \put( -6,98){$A_1$}   \put( 15,40){$A_2$}
   \put( 99,65){$A_3$}   \put( 87,11){$A_4$}
   \put(134,28){$A_5$}   \put(212,25){$A_6$}

   \put(124,190){$B_1$}   \put(104,138){$B_2$}
   \put(198,185){$B_3$}   \put(146,112){$B_4$}
   \put(232,170){$B_5$}   \put(234,117){$B_6$}

   \put( 64,144){$l_1$}   \put( 68, 97){$l_2$}
   \put(119, 92){$l_3$}   \put(107, 41){$l_4$}
   \put(186, 65){$l_5$}   \put(216, 58){$l_6$}

  \end{picture}
$$
\caption{A Stewart Platform\label{fig:stewart}}
\end{figure}

The position of a rigid body in ${\mathbb R}^3$ has 6 degrees of freedom.
This is exploited in robotics, giving rise to the Stewart-Gough
platform~\cite{Go57,St65}.
Specifically, suppose we have 6 fixed points $A_1,A_2,\ldots, A_6$ in space
and 6 points $B_1,B_2,\ldots, B_6$ on a rigid body $B$ (the framework of
Figure~\ref{fig:stewart}).
The body is controlled by varying each distance $l_i$ between the
fixed point $A_i$ and the point $B_i$ on $B$.
This may be accomplished by attaching rigid actuators to
spherical joints located at the points $B_i$, or by suspending the platform from
cables.

Given a position of the body $B$ in space, the distances $l_1,l_2,\ldots, l_6$ are
uniquely determined.
A fundamental problem is the inverse question:
Given a platform (positions of the $A_i$ fixed and the relative positions of
the $B_i$ specified) and a sextuple of distances $l_1,l_2,\ldots,l_6$, what is the
position of the platform in space?

It had long been understood that several positions were possible for a given
sextuple of lengths.
An early work in 1897 showed there could be as many as 16 different
positions~\cite{Br1897}. 
This leads to the following enumerative problem.

\begin{ques}
 For a given Stewart platform, how many (complex) positions are there for a
 generic choice of the distances $l_1,l_2,\ldots,l_6$?
 How many of these can be real?
\end{ques}

In the early 1990's, several approaches (numerical
experimentation~\cite{Ra91}, intersection theory~\cite{RV95}, Gr\"obner
bases~\cite{La93}, resultants~\cite{Mo93}, and algebra~\cite{Mo94}) each
showed that there are 40 complex positions of a general Stewart platform.
The obviously practical question of how many positions could be real 
remained open until 1998, when Dietmaier introduced a novel method to
find a value of the distances $l_1,l_2,\ldots,l_6$ with all 40 positions real.

\begin{thm}[Dietmaier~\cite{Di98}]
  All $40$ positions can be real.
\end{thm}

Dietmaier's method will find future applications to other problems of this kind.
He began with a formulation of the problem as a system of equations depending
upon the distances $l_1,l_2,\ldots,l_6$.
An initial solution for a given instance of the distances
gave 6 real solutions and 17 pairs of complex conjugate solutions.
He then used an ingenious algorithm to vary the distances in search of a
configuration with all 40 solutions real.

\begin{figure}[htb]
$$
  \begin{picture}(170,100)
   \put( 9,0){\epsfxsize=2.2in\epsfbox{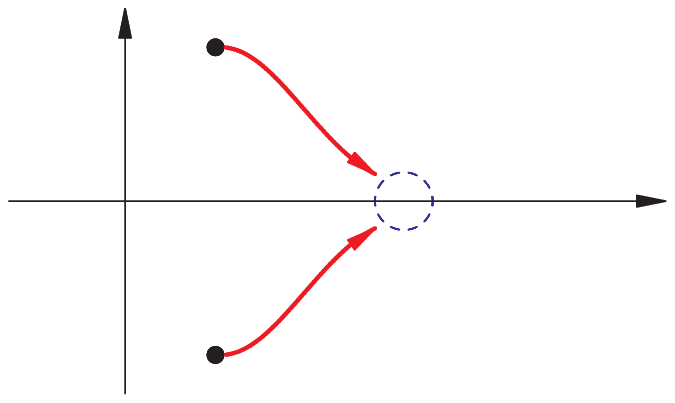}}
   \put( 0,75){Im$(x)$}   \put(125,60){Re$(x)$}
   \put(44,8){$\overline{x}$}
   \put(44,85){$x$}
  \end{picture}
  \qquad\begin{picture}(15,100)
          \put(0,47){$\Longrightarrow$}
        \end{picture}\qquad
  \begin{picture}(170,100)
   \put( 9,0){\epsfxsize=2.2in\epsfbox{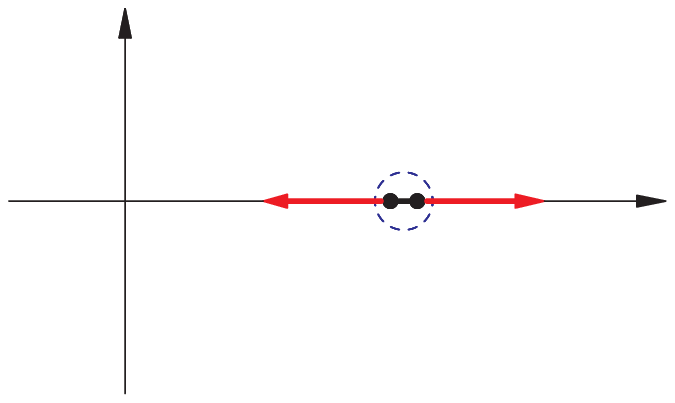}}
   \put( 0,75){Im$(x)$}   \put(125,60){Re$(x)$}
   \put( 79,40){$x_1$}
   \put(118,40){$x_2$}
  \end{picture}
$$
\caption{The two stages of Dietmaier's algorithm.\label{fig:Dietmaier}}
\end{figure}

This algorithm systematically varies the distances with the intention of
increasing the number of real solutions.
It proceeds in two stages.
In the first stage, a pair of complex conjugate solutions ($x,\overline{x}$) are
driven together, eventually creating a double solution,
while at the same time the existing real solutions are kept bounded away from
one another.
At the formation of a double (necessarily real) solution, the distances are
further incremented to create two new nearby real solutions ($x_1, x_2$),
which are then driven apart in the second stage.
This procedure is repeated again with another pair of complex conjugate
solutions, and {\it et cetera}.
Figure~\ref{fig:Dietmaier} ilustrates the two stages.

In each stage, Dietmaier accomplishes the given task (eg.~colliding conjugate
solutions) by linearizing the system at the current solutions and then  
solving a linear program for the optimal increment of the distances for
the given goal. 
Changing the distances, he uses Newton's method beginning with the
current solutions to find soutions for the new set of distances, and then
repeats this procedure until the goal is acheived (eg.~the conjugate pair
collides). 
This is an application of numerical homotopy continuation~\cite{Ver99}. 

While there is no guarantee that this method will even successfully collide
two conjugate solutions, Dietmaier uses it to find a sextuple of distances
with all 40 solutions real.
While at each step the solutions are only numerical
approximations to the actual solutions, the condition number $N$
guarantees the existence of a genuine solution within $1/N$ of each
approximate solution. 
Since the approximate real solutions were separated by more than 
$2/N$, the requirement that non-real solutions occur in complex
conjugate pairs forced these genuine solutions to be real.

\subsection{Real rational cubics through 8 points in 
            ${\mathbb P}^2_{\mathbb R}$}\label{sec:ratcubics}
There are 12 singular (rational) cubic curves containing 8 general points in the
plane.

\begin{thm}[Degtyarev and 
     Kharlamov~{\cite[Proposition 4.7.3]{DeKh00}}]\label{thm:cubics}
 Given $8$ general points in ${\mathbb P}^2_{\mathbb R}$, there are at least $8$
 singular real cubics containing them, and there are choices of the $8$ points
 for which all\/ $12$ singular cubics are real.
\end{thm}

\begin{proof}
 Since a cubic equation in the plane has 10 coefficients, 
 the space of cubics is identified with ${\mathbb P}^9$.
 The condition for a plane cubic to contain a given point is linear in these
 coefficients. 
 Given 8 general points, these linear equations are independent and so there is
 a pencil (${\mathbb P}^1$) of cubics containing 8 general points in 
 ${\mathbb P}^2$.

 Two cubics $P$ and $Q$ in this pencil meet
 transversally in 9 points.
 Since curves in the pencil are given by $aP+bQ$ for $[a,b]\in{\mathbb P}^1$,
 any two curves in the pencil meet transversally in these 9 points.
 Let $Z$ be ${\mathbb P}^2$ blown up at these same 9 points.
 We have a map
$$
  f\ \colon\ Z\ \longrightarrow\ {\mathbb P}^1\,,
$$
 where $f^{-1}([a,b])$ is the cubic curve defined by the polynomial $aP+bQ$.

 Consider the Euler characteristic $\chi(Z)$ of $Z$ first over ${\mathbb C}$
 and then over ${\mathbb R}$.
 Blowing up a smooth point on a surface replaces it with a 
 ${\mathbb P}^1_{\mathbb C}$ and thus increases the Euler characteristic by 1.
 Since $\chi({\mathbb P}^2_{\mathbb C})=3$, we see that $\chi(Z)=12$.
 The general fiber of $f$ is a smooth plane cubic which is homeomorphic to
 $S^1\times S^1$, and so has Euler characteristric 0.
 Thus only the singular fibers of $f$ contribute to the Euler characteristic
 of $Z$. 
 Assume that the 8 points are in general position so there are only nodal 
 cubics in the pencil.
 A nodal cubic has Euler characteristic 1.
 Thus there are 12 singular fibers of $f$ and hence 12 singular cubics meeting 8
 general points in ${\mathbb P}^2_{\mathbb C}$.\smallskip

 Consider now the Euler characteristic of $Z_{\mathbb R}$.
 Blowing up a smooth point on a real surface replaces the point by 
 ${\mathbb P}^1_{\mathbb R}=S^1$, and hence decreases the Euler characterstic
 by 1.
 Since $\chi({\mathbb P}^2_{\mathbb R})=1$, we have 
 $\chi(Z_{\mathbb R})=1-9=-8$.
 A nonsingular real cubic is homeomorphic to either one or two disjoint copies
 of $S^1$, and hence has Euler characteristic 0.
 Again the Euler characteristic of $Z_{\mathbb R}$ is carried by its singular
 fibers.
 There are two nodal real cubics; either the node has two real
 branches or two complex conjugate branches so that the singular point is
 isolated. 
 Call these curves real nodal and complex nodal, respectively.
 They are displayed in Figure~\ref{fig:nodes}.
\begin{figure}[htb]
  $$
    \epsfxsize=1.35in \epsfbox{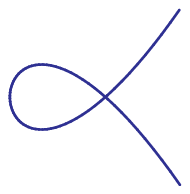}\qquad \qquad 
    \epsfxsize=1.35in \epsfbox{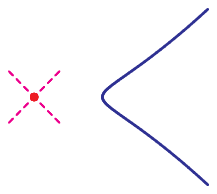}
  $$
  \caption{A real nodal and complex nodal curve\label{fig:nodes}}
\end{figure}
 The real nodal curve is homoemorphic to a figure 8 and has Euler
 characteristic $-1$, while the complex nodal curve is the union of a $S^1$
 with a point and so has Euler characteristic $1$.

 Among the singular fibers, we have
 \begin{eqnarray*}
  -8&=& \#\mbox{complex nodal}\; -\; \#\mbox{real nodal}\,,\quad\mbox{with}\\
  12&\geq& \#\mbox{complex nodal}\; +\; \#\mbox{real nodal}\,.
 \end{eqnarray*}
 Thus there are at least 8 real nodal curves containing 8 general points in
 ${\mathbb P}^2_{\mathbb R}$.
 The pencil of cubics containing the 2 complex nodal cubics of
 Figure~\ref{fig:9-int} has 10 real nodal cubics.
 Thus there are 12 real rational cubics containing any 8 of the 12
 points in Figure~\ref{fig:9-int}.
\begin{figure}[htb]
 \begin{minipage}[c]{3.in}
   \begin{eqnarray*}
    \frac{(y+9x-28)^2}{4}&=&(x-1)x^2\\
    \frac{(x+10y-28)^2}{4}&=&y(y+1)^2\rule{0pt}{23pt}
   \end{eqnarray*}
 \end{minipage}\ 
 \begin{minipage}[c]{2.2in}
 $$
   \epsfxsize=2in \epsfbox{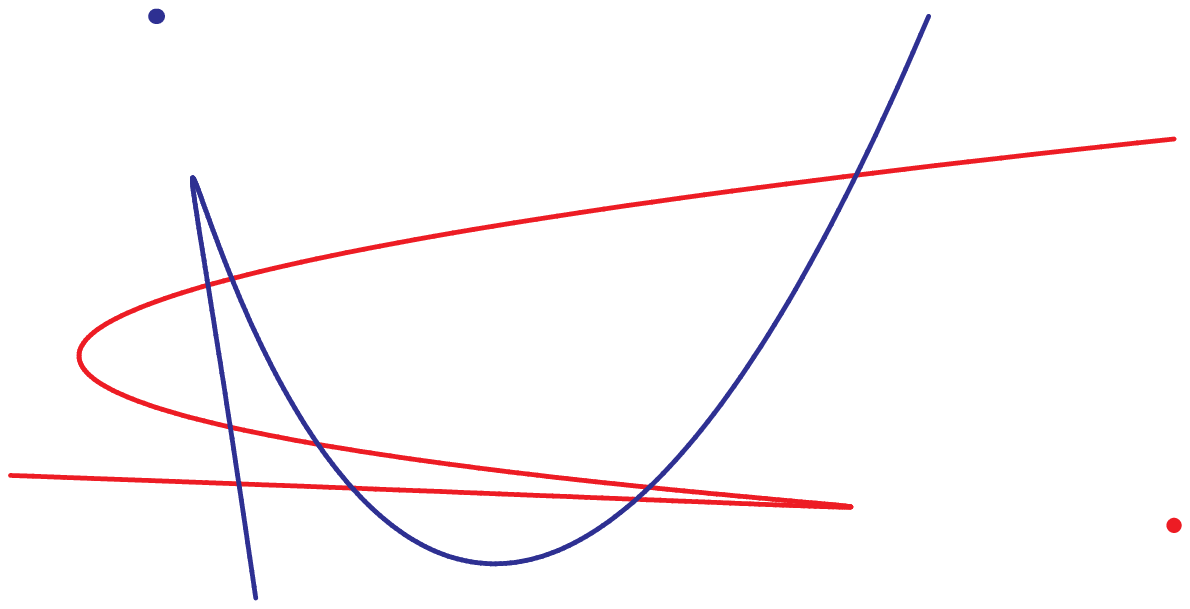}
 $$
 \end{minipage}
\caption{Complex nodal curves meeting in 9 points\label{fig:9-int}}
\end{figure}

\end{proof}

\begin{remark}
 This classical problem of 12 plane cubics containing 8 points generalizes to
 the problem of enumerating rational plane curves of degree $d$ containing
 $3d-1$ points.
 Let $N_d$ be the number of such curves, which satisfies the
 recursion~\cite{FP97}
$$
  N_d\ =\ \sum_{d_1+d_2=d}  N_{d_1}N_{d_2}\left[
              d_1^2d_2^2\binom{3d-4}{3d_1-2}
             -d_1^3d_2\binom{3d-4}{3d_1-1}\right]\,.
$$
 The values $N_1=N_2=1$ are trivially fully real, and we have
 just seen that $N_3=12$ is fully real.
 The next case of $N_4=620$ (computed by Zeuthen~\cite{Ze1873}) 
 seems quite challenging.
\end{remark}

\begin{remark}
 The most 
 interesting feature of Theorem~\ref{thm:cubics} is the existence of a lower
 bound on the number of real solutions, which is a new phenomenon.
 In Section~\ref{sec:lower} we shall see evidence that this may be a pervasive
 feature of this field. 
\end{remark}

\subsection{Common tangent lines to $2n-2$ spheres in 
 ${\mathbb R}^n$}\label{sec:tangent} 
Consider the following.

\begin{ques}
 How many common tangent lines are there to $2n-2$ spheres in ${\mathbb R}^n$?
\end{ques}

For example, when $n=3$, how many common tangent lines are there to four
spheres in ${\mathbb R}^3$?
Despite its simplicity, this question does not seem to have been asked
classically, but rather arose in discrete and computational
geometry\frankfootnote{The question of the maximal number of (real) common tangents
to 4 balls was first formulated by David Larman~\cite{La90} at DIMACS in 1990.}.
The case $n=3$ was solved by Macdonald, Pach, and Theobald~\cite{MPT01} and
the general case more recently\frankfootnote{This was solved during the DIMACS
Workshop on Algorithmic and Quantitative Aspects of Real Algebraic Geometry in
Mathematics and Computer Science.}~\cite{STh01}. 

\begin{thm}\label{thm:12lines}
 $2n-2$ general spheres in ${\mathbb R}^n$ ($n\geq 3$) have $3\cdot 2^{n-1}$
 complex common tangent lines, and there are $2n-2$ such spheres with all
 common tangent lines real.
\end{thm}

\begin{figure}[htb]
$$
  \epsfysize=2.5in \epsfbox{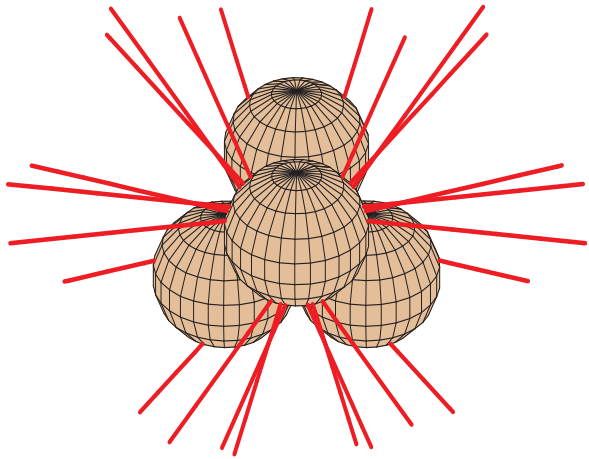}
$$
\caption{Four spheres with 12 common tangents\label{fig:12lines}}
\end{figure}

Figure~\ref{fig:12lines} displays a configuration of 4 spheres in
${\mathbb R}^3$ with 12 real common tangents.
The number $2n-2$ is the dimension of the space
of lines in ${\mathbb R}^n$ and is necessary for there to be finitely many
common tangents.

Represent a line in ${\mathbb R}^n$ by a point $p$ on the line and its
direction vector $v\in{\mathbb P}^{n-1}$.
Imposing the condition
 \begin{equation}\label{eq:dot}
    p\cdot v\ =\ 0\,,
 \end{equation}
makes this representation unique.
Here $\cdot$ is the usual Euclidean dot product.
Write $v^2$ for $v\cdot v$.
A line $(p,v)$ is tangent to the sphere with center $c$ and radius $r$ when
$$
  v^2p^2-2v^2p\cdot c + v^2c^2 - (v\cdot c)^2 - r^2v^2\ =\ 0\,.
$$

Suppose we have $2n-2$ spheres with centers $c_1,c_2,\ldots,c_{2n-2}$ and
corresponding radii $r_1,r_2,\ldots,r_{2n-2}$.
Without any loss of generality, we may assume that the last sphere is centered
at the origin and has radius $r$.
Then its equation is
 \begin{equation}\label{eq:simple}
  v^2p^2-v^2r^2\ =\ 0\,.
 \end{equation}
Subtracting this from the equations for the other
spheres, we obtain the equations
 \begin{equation}\label{eq:rest}
  2v^2p\cdot c_i\ =\ v^2c_i^2-(v\cdot c_i)^2 - v^2(r_i^2-r^2)
 \end{equation}
for $i=1,2,\ldots,2n-3$.

These last equations are linear in $p$.
If the centers $c_1,c_2,\ldots,c_n$ are linearly independent 
(which they are, by
our assumption on generality), then we use the
corresponding equations to solve for $v^2p$ 
as a homogeneous quadratic in $v$.
Substituting this into the equations~(\ref{eq:dot}) and~(\ref{eq:simple})
gives a cubic and a quartic in $v$, and substituting the expression for $v^2p$
into~(\ref{eq:rest}) for $i=n+1,n+2,\dots,2n{-}3$ gives $n{-}3$ quadratics
in $v$.
By B\'ezout's Theorem, if there are finitely many complex solutions to these
equations, their number is bounded by
$ 3\cdot 4\cdot 2^{n-3} = 3\cdot 2^{n-1}$.

This upper bound is attained with all solutions real.
Suppose that the spheres all have the same radius, $r$, and the first four
have centers 
 \begin{eqnarray*}
  c_1 &:=& (\hh 1,\hh 1,\hh 1,\ 0,\ldots,0)\,,\\
  c_2 &:=& (\hh 1,-1,-1,\ 0,\ldots,0)\,,\\
  c_3 &:=& (-1,\hh 1,-1,\ 0,\ldots,0)\,,\\
  c_4 &:=& (-1,-1,\hh 1,\ 0,\ldots,0)\,,
 \end{eqnarray*}
and subsequent centers are at the points $\pm ae_j$  for $j=4,5,\ldots,n$,
where $e_1,e_2,\ldots,e_n$ is the standard basis for ${\mathbb R}^n$.
Let $\gamma:=a^2(n-1)/(a^2+n-3)$, which is positive.

\begin{thm}[{\cite[Theorem 5]{STh01}}]\label{thm:real-cmplx}
 When 
$$
  a\,r\,(r^2-3)\,(3-\gamma)\,(a^2-2)\,(r^2-\gamma)\, 
  \left((3-\gamma)^2 +4\gamma - 4r^2\right)\ \neq \ 0\,,
$$
 there are exactly $3\cdot 2^{n-1}$ complex lines tangent to the spheres.
 If we have
\begin{enumerate}
\item[(a)] ${\displaystyle 
         \frac{1}{4}\left(3-\gamma\right)^2+\gamma\ >\  
                r^2\ >\ \gamma}$\ and
\item[(b)] ${\displaystyle \frac{n-1}{n-4} + 2\ >\ a^2\ >\ 2}$\,, 
           \rule{0pt}{20pt} 
\end{enumerate}
then all the $3\cdot 2^{n-1}$ lines are in fact real.
\end{thm}

Theorem~\ref{thm:12lines} is false when $n=2$.  
There are 4 lines tangent to 2 circles in the plane, and all
may be real.
The argument given for Theorem~\ref{thm:12lines} fails because 
the centers of the spheres do not affinely span 
${\mathbb R}^2$.
This case of $n=2$ does generalize, though.

\begin{thm}[Megyesi~\cite{Me01}]
    Four unit spheres in ${\mathbb R}^3$ whose centers are coplanar but
    otherwise general have $12$ common complex tangents.
    At most $8$ of these $12$ are real.
\end{thm}

\begin{remark}
 This problem of common tangents to 4 spheres with equal radii
 and coplanar centers gives an example of an enumerative geometric
 problem that is not fully real.
 We do not feel this contradicts the observation that there are no enumerative
 problems not known to be fully real, as the spheres are not sufficiently
 general. 
\end{remark}

\section{Schubert Calculus}\label{sec:SchubertCalculus}

The classical Schubert calculus of enumerative geometry is concerned with
questions of enumerating linear subspaces of a vector space or projective
space that satisfy incidence conditions imposed by other linear subspaces.
A non-trivial instance is the question posed at the
beginning of Section~3.

\begin{ques}
 How many lines in space meet four general lines $\ell_1,\ell_2,\ell_3$, and
 $\ell_4$?
\end{ques}

\begin{figure}[htb]
$$
  \setlength{\unitlength}{0.9pt}
  \setlength{\unitlength}{1pt}
  \begin{picture}(270,160)(-5,-5)
   \put(  5,  5){\epsfysize=2in \epsfbox{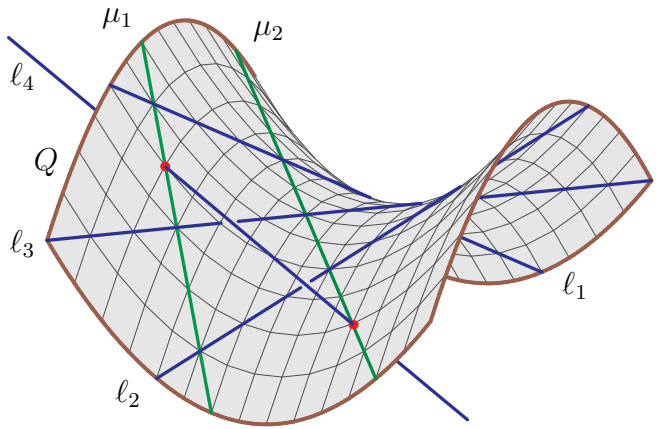}}
   \put( 35,  8){$\ell_2$}
   \put( -5, 65){$\ell_3$}
   \put(204, 49){$\ell_1$}
   \put( -5,128){$\ell_4$}
   \put( 30,152){$\mu_1$}
   \put( 87,147){$\mu_2$}
   \put(  4,96){$Q$}
  \end{picture}
$$
\caption{The two lines meeting four general lines in space.\label{fig:4lines}}
\end{figure}

Three pairwise skew lines  $\ell_1,\ell_2$, and $\ell_3$ lie on a unique
smooth quadric surface $Q$.
There are two families of lines that foliate $Q$---one family includes 
$\ell_1,\ell_2$, and $\ell_3$ and the other consists of the lines meeting each
of $\ell_1,\ell_2$, and $\ell_3$.
The fourth line $\ell_4$ meets $Q$ in two points, and each of these points
determines a line in the second family.
Thus there are 2 lines $\mu_1,\mu_2$ in space that meet general lines 
$\ell_1,\ell_2,\ell_3$, and $\ell_4$.
Figure~\ref{fig:4lines} shows this configuration.
Figure~\ref{fig:4lines} also shows how the two lines can be real---if $\ell_4$
meets $Q$ in two real points.
(The two lines are complex when $\ell_4$ meets $Q$ in two complex conjugate
points.) 

The classical Schubert calculus is a vast generalization of this problem of
four lines.
In the 1980's Robert Speiser suggested to Fulton that the classical
Schubert calculus may be a good testing ground for
Question~\ref{ques:reality}.
This was also considered by Chiavacci and Escamilla-Castillo~\cite{CEC88}.
We will discuss increasingly more general versions of the Schubert calculus,
and the status of Question~\ref{ques:reality} for each.

Consider first more general problems involving lines.
The space of lines in ${\mathbb P}^n$ is a smooth projective variety of
dimension $2n{-}2$ called the Grassmannian of lines in ${\mathbb P}^n$.
The set of lines meeting a linear subspace $L$ of dimension $n{-}1{-}l$ has
codimension $l$ in the Grassmannian.
Thus given general linear subspaces $L_1,L_2,\ldots,L_s$ of ${\mathbb P}^n$  with
$\dim L_i=n-1-l_i$ where $l_1+l_2+\cdots+l_s=2n-2$, we expect (and there are
indeed) finitely many lines in ${\mathbb P}^n$ meeting each linear subspace 
$L_1,L_2,\ldots,L_s$.
Schubert~\cite{Sch1886a} discovered algorithms for computing this number
$d(l_1,l_2,\ldots,l_s)$ of lines.
For example, if each $l_i=1$, so that $s=2n-2$, then this number is the
$n$th Catalan number\frankfootnote{This indexing of the Catalan numbers is shifted
from that of some other authors.}
 \begin{equation}\label{eq:Catalan}
  C_n\ =\ \frac{1}{n}\binom{2n-2}{n-1}\,.
 \end{equation}
%

Enumerative problems of lines in ${\mathbb P}^n$ meeting 
general linear subspaces furnished the first infinite family of non-trivial
enumerative problems known to be fully real.

\begin{thm}[{\cite[Theorem C]{So97a}}]\label{thm:reality}
 Given positive integers $l_1,l_2,\ldots,l_s$ with $l_1+l_2+\cdots+l_s=2n-2$,
 there  exist linear subspaces $L_1,L_2,\ldots,L_s$ of\/ 
 ${\mathbb R}{\mathbb P}^n$ 
 with $\dim L_i=n+1-l_i$ such that there are exactly $d(l_1,l_2,\ldots,l_s)$ 
 complex lines meeting each subspace $L_i$, and each of these lines are real.
\end{thm}

\subsection{The special Schubert calculus}\label{sec:SpSchCalc}
More generally, we may ask how many linear
subspaces of a fixed dimension meet general linear subspaces.
We formulate this question in terms of linear subspaces of a vector
space.\smallskip 

The set of $k$-dimensional subspaces ($k$-planes) of an $n$-dimensional vector
space forms the {\it Grassmannian of $k$-planes in $n$-space}, $\Gr(k,n)$, 
a smooth projective variety of dimension $k(n{-}k)$.
Those $k$-planes meeting a linear subspace $L$ of dimension $n{-}k{+}1{-}l$ 
non-trivially (that is, the intersection has positive dimension) form the special
Schubert subvariety $\Omega(L)$ of $\Gr(k,n)$ which has 
codimension $l$.
The special Schubert calculus is concerned with the following question.

\begin{ques}\label{ques:spSchub}
 Given general linear subspaces $L_1,L_2,\ldots,L_s$ of ${\mathbb C}^n$ with 
 $\dim L_i=n-k+1-l_i$ where $l_1+l_2+\cdots+l_s=k(n-k)$, how many $k$-planes $K$
 meet each subspace $L_i$ non-trivially, that is, satisfy
 \begin{equation}\label{eq:special}
   K\cap L_i\ \neq\ \{0\}\qquad i=1,2,\ldots,s\,?
 \end{equation}
\end{ques}

The condition~(\ref{eq:special}) is expressed in the global geometry of
$\Gr(k,n)$ as the number of points in the intersection of the special Schubert
varieties
 \begin{equation}\label{eq:special-int}
   \Omega(L_1)\cap \Omega(L_2)\cap\cdots\cap\Omega(L_s)\,,
 \end{equation}
when the intersection is transverse.
(A general theorem of Kleiman~\cite{MR50:13063} guarantees transversality when
the $L_i$ are in general position, and also implies transversality for the
other intersections considered in this section.)

There are algorithms due to Schubert~\cite{Sch1886b} (when each $l_i=1$) and
Pieri~\cite{Pi1893} to compute the expected number of solutions.
When each $l_i=1$, Schubert~\cite{Sch1886c} showed that the number of
solutions is equal to 
$$
 d(n,k)\ :=\ \frac{1!\,2!\,\cdots\,(k\!-\!1)!\cdot[k\,(n\!-\!k)]!}
  {(n\!-\!k)!\,(n\!-\!k\!+\!1)!\cdots (n\!-1)!}\ . 
$$

A line in ${\mathbb P}^n$ is a 2-plane in $(n+1)$-space
and two linear subspaces in ${\mathbb P}^n$ meet if and only if the
corresponding linear subspaces in $(n+1)$-space have a non-trivial 
intersection. 
Thus the problem of lines in projective space corresponds to the case $k=2$ of
the special Schubert calculus.
While the geometric problem generalizes easily from $k=2$ to arbitrary values
of $k$, the proof of Theorem~\ref{thm:reality} does not.
There is, however, a relatively simple argument that this special Schubert
calculus is fully real.

\begin{thm}[{\cite[Theorem 1]{So99a}}]\label{thm:special}
 Suppose $n>k>0$ and $l_1,l_2,\ldots,l_s$ are positive integers with
 $l_1+l_2+\cdots+l_s=k(n-k)$. 
 Then there are linear subspaces $L_1,L_2,\ldots,L_s$ of\/ ${\mathbb R}^n$ in
 general position with $\dim L_i=n-k+1-l_i$ such that
 each of the a priori complex $k$-planes $K$ satisfying~$(\ref{eq:special})$
 are in fact real.
\end{thm}

We present an elementary proof of this result in the important special case
when each $l_i=1$ so that the conditions are simple, meaning each
$\Omega(L_i)$ has codimension 1.
This proof generalizes to show that some other classes of enumerative problems
in the Schubert calculus are fully real (see Sections~\ref{sec:quantum}
and~\ref{sec:flags}). 
This generalization constructs sufficiently many real solutions using a
limiting argument, as in Section~\ref{sec:sparse}.
Just as the arguments of Section~\ref{sec:sparse} were linked to the homotopy
algorithms of Huber and Sturmfels, the proof of Theorem~\ref{thm:special}
leads to numerical homotopy methods for solving these
problems~\cite{HSS98,HV00}. 

We develop further geometric properties of Grassmann varieties.
The $k$th exterior power of the embedding
$K\ \to\ {\mathbb C}^n$
of a $k$-plane $K$ into ${\mathbb C}^n$ gives the embedding
 \begin{equation}\label{eq:exter-prod}
  {\mathbb C}\ \simeq\ \wedge^k K\ \longrightarrow\ \wedge^k\,{\mathbb C}^n\,,
 \end{equation}
whose image is a line in $\wedge^k{\mathbb C}^n$ and thus a point in the
projective space 
${\mathbb P}(\wedge^k{\mathbb C}^n)\simeq{\mathbb P}^{\binom{n}{k}-1}$.
This point determines the $k$-plane $K$ uniquely.
The {\it Pl\"ucker embedding} is the resulting projective embedding of the
Grassmannian
$$
  \Gr(k,n)\ \longrightarrow\ {\mathbb P}^{\binom{n}{k}-1}\,.
$$

The $\binom{n}{k}$ homogeneous Pl\"ucker coordinates for the
Grassmannian in this embedding are realized concretely as follows.
Represent a $k$-plane $K$ as the row space of a $k\times n$ matrix, also
written $K$.
A maximal minor of $K=(a_{ij})$ is the determinant of a $k\times k$ submatrix
of $K$:
Given a choice of columns 
$\alpha\colon 1\leq \alpha_1<\alpha_2<\cdots<\alpha_k\leq n$, set
$$
   p_\alpha(K)\ :=\ \det\left[\begin{array}{cccc}
           a_{1,\alpha_1}&a_{1,\alpha_2}&\cdots& a_{1,\alpha_k}\\
           a_{2,\alpha_1}&a_{2,\alpha_2}&\cdots& a_{2,\alpha_k}\\
                \vdots   &     \vdots   &\ddots&    \vdots\\
           a_{k,\alpha_1}&a_{k,\alpha_2}&\cdots& a_{k,\alpha_k}
                \end{array}\right]\ .
$$
The vector $(p_\alpha(K))$ of $\binom{n}{k}$ maximal minors of $K$ defines the
map~(\ref{eq:exter-prod}) giving {\it Pl\"ucker coordinates} for $K$.
Let $\binom{[n]}{k}$ be the collection of these indices of Pl\"ucker
coordinates. 

The indices $\binom{[n]}{k}$ have a natural Bruhat order
$$
  \beta\leq \alpha\ \Longleftrightarrow\ 
     \beta_j\leq\alpha_j\quad\mbox{for}\quad j=1,2,\ldots,k\,.
$$
The Schubert variety $\Omega_\alpha$  is 
 \begin{equation}\label{eq:sch-def}
  \Omega_\alpha\ =\ \{ K\in\Gr(k,n)\ \mid\ 
    p_\beta(K)=0\ \mbox{for}\ \beta\not\leq\alpha\}\,.
 \end{equation}
This has dimension $|\alpha|:=\sum_j(\alpha_j-j)$.

The relevance of the Pl\"ucker embedding to Question~\ref{ques:spSchub} when
$l_i=1$ is seen as follows.
Let $L$ be a $(n{-}k)$-plane, represented as the row space of a $(n{-}k)$ by
$n$ matrix, also written $L$.
Then a general $k$-plane $K$ meets $L$ non-trivially if and only if
$$
  \det\left[\begin{array}{c}K\\L\end{array}\right]\ =\ 0\,.
$$
Laplace expansion along the rows of $K$ gives
 \begin{equation}\label{eq:laplace}
  0\ =\ \det\left[\begin{array}{c}K\\L\end{array}\right]
   \ =\ \sum_{\alpha\in\binom{[n]}{k}} p_\alpha(K) L_\alpha\,,
 \end{equation}
where $L_\alpha$ is the
appropriately signed minor of $L$ given by the columns complementary to
$\alpha$.
Hence the set $\Omega(L)$ of $k$-planes that meet the $(n{-}k)$-plane $L$
non-trivially is a hyperplane section of the Grassmannian in its Pl\"ucker
embedding.

Thus the set of $k$-planes meeting $k(n{-}k)$ general
$(n{-}k)$-planes non-trivially is a complementary linear section
of the Grassmannian, and so the number $d(k,n)$ of such $k$-planes is the
degree of the Grassmannian in its Pl\"ucker embedding.
More generally, if $\alpha\in\binom{[n]}{k}$ and $L_1,L_2,\ldots,L_{|\alpha|}$
are general $(n{-}k)$-planes, then the number of points in the intersection 
 \begin{equation}\label{eq:schub-int}
  \Omega_\alpha\,\cap\,\Omega(L_1)\,\cap\,\Omega(L_2)\,
    \cap\,\cdots\,\cap\,\Omega(L_{|\alpha|})
 \end{equation}
is the degree $d(\alpha)$ of the Schubert variety $\Omega_\alpha$, which we
now compute.

An intersection $X\cap Y$ is {\it generically transverse} if $X$ and $Y$ meet
transversally along an open subset of every component of $X\cap Y$.
When $\beta,\alpha\in\binom{[n]}{k}$ satisfy $\beta<\alpha$ but there
is no index $\gamma$ with $\beta<\gamma<\alpha$, then we write
$\beta\lessdot\alpha$. 
The following fact is elementary and due to Schubert.

\begin{thm}\label{thm:Schubert-cover}
 Let $\alpha\in\binom{[n]}{k}$ and set $H_\alpha$ to be the hyperplane defined
 by $p_\alpha=0$.
 Then
$$
  \Omega_\alpha \cap H_\alpha\ =\ \bigcup_{\beta\lessdot\alpha}\Omega_\beta\,,
$$
 and the intersection is generically transverse.
\end{thm}

\begin{figure}[htb]
$$
  \epsfysize=3.1in \epsfbox{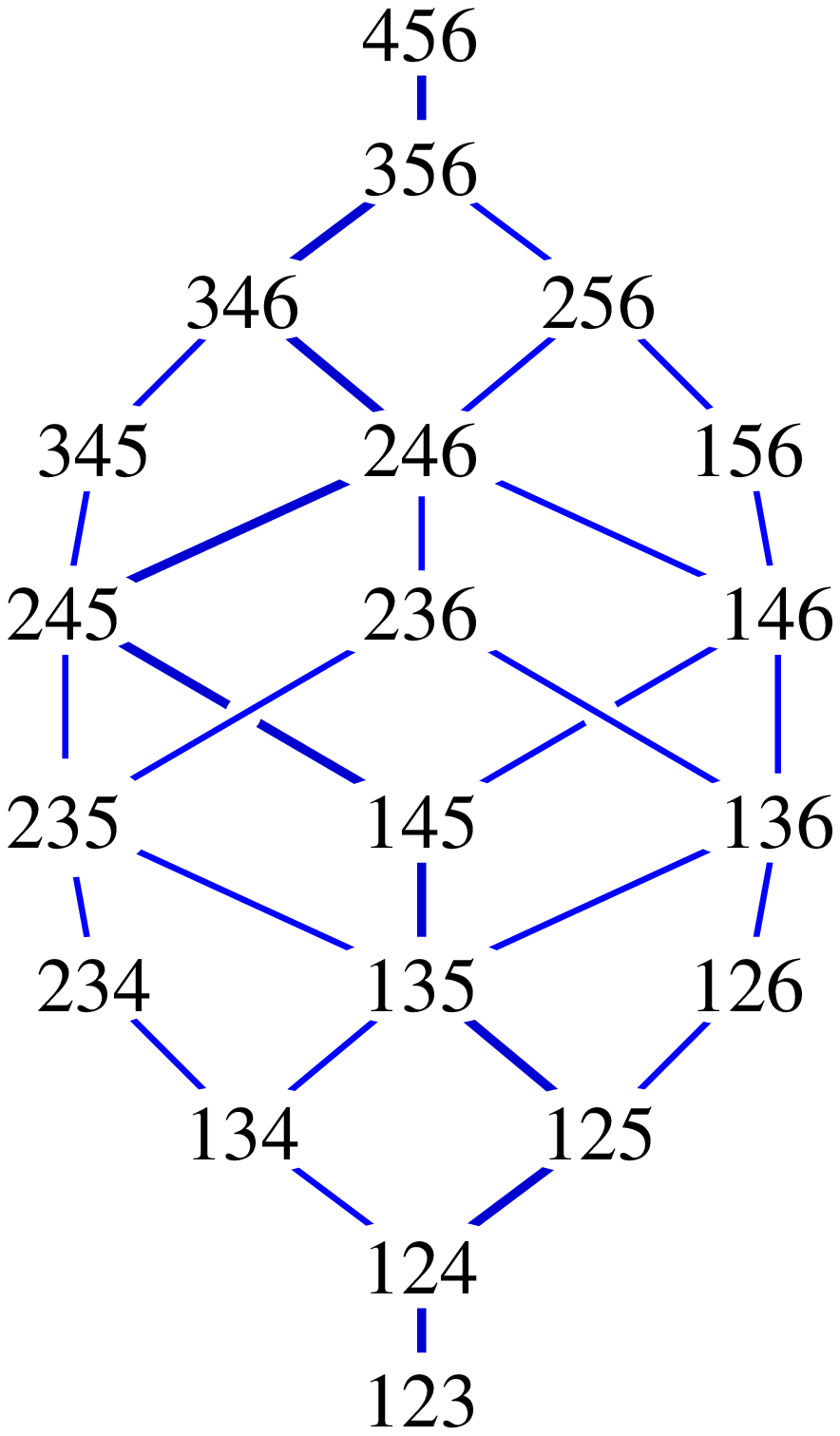}\qquad
  \epsfysize=3.1in \epsfbox{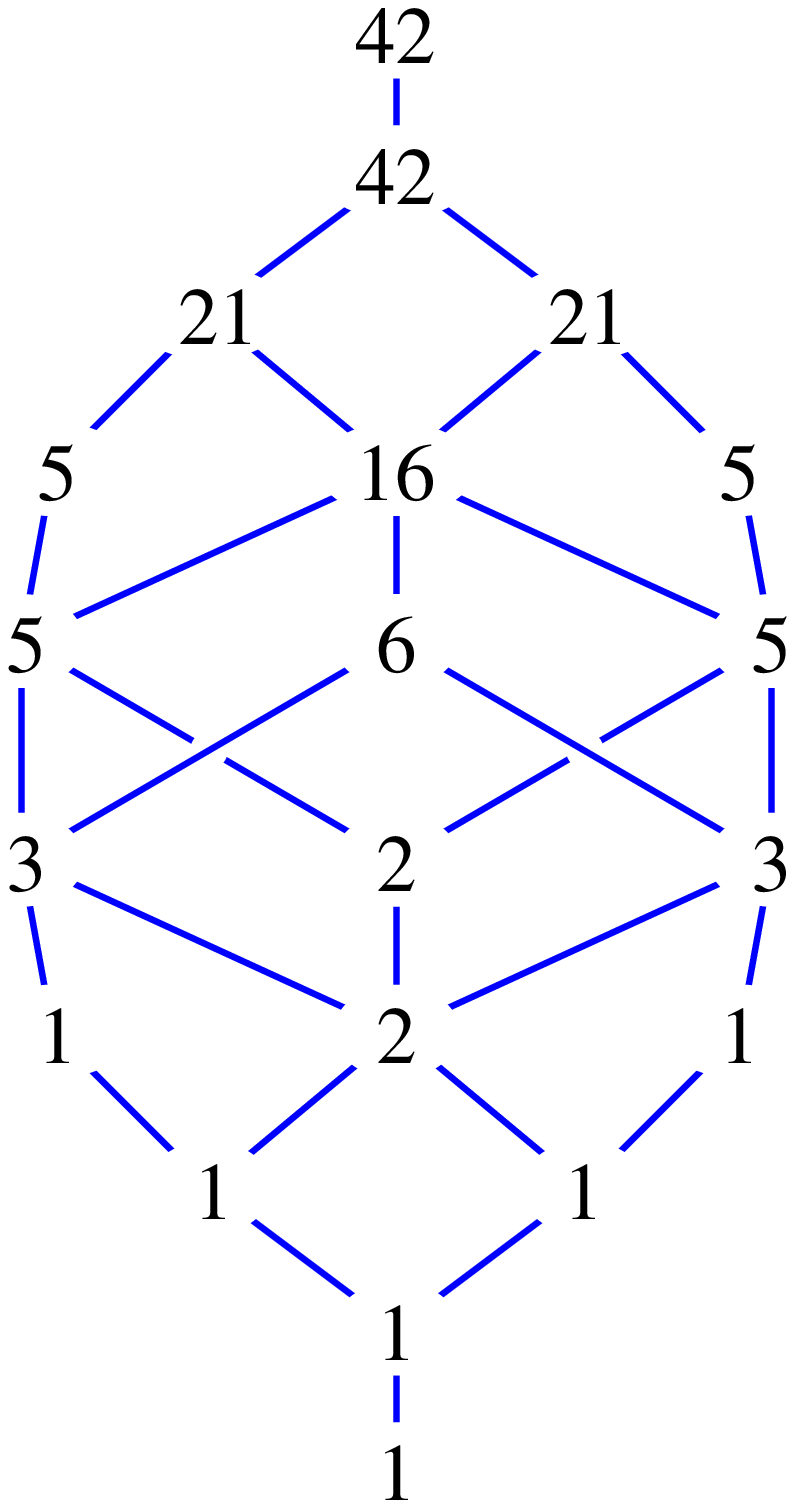}
$$
\caption{Bruhat order and degrees of Schubert varieties, $k=3$, $n=6$}
\label{fig:BrDeg}
\end{figure}

In fact this intersection is transverse along 
$\Omega^\circ_\beta:=\Omega_\beta-\bigcup_{\delta<\beta}\Omega_\delta$.
We obtain the recursion for the degree $d(\alpha)$ of the Schubert
variety $\Omega_\alpha$
$$
   d(\alpha)\ =\ \sum_{\beta\lessdot\alpha} d(\beta)\,.
$$
Since the minimal Schubert variety is a point (which has degree 1), this gives 
a conceptual formula for $d(\alpha)$.  
Let $\hat{0}=(1,2,\ldots,k)$ be the minimal element in the Bruhat order.
 \begin{equation}\label{eq:deg-BO}
  d(\alpha)\ =\ \mbox{the number of saturated chains in the Bruhat order from
  $\hat{0}$ to $\alpha$}\,.
 \end{equation}
Figure~\ref{fig:BrDeg} displays both the
Bruhat order for $k=3$ and $n=6$ (on the left) and the degrees of the
corresponding Schubert varieties (on the right).

Consider the action of the non-zero real numbers ${\mathbb R}^\times$
on ${\mathbb R}^n$
 \begin{equation}\label{eq:action}
   t.e_j\ :=\ t^j\cdot e_j\,,
 \end{equation}
where $t\in{\mathbb R}^\times$ and $e_1,e_2,\ldots,e_n$ is a basis for
${\mathbb R}^n$ (corresponding to the rows of the $n\times n$ identity matrix).
Let $L$ be a $(n-k)$-plane.
By~(\ref{eq:laplace}), the equation for a $k$-plane $K$ to meet
$t.L$ non-trivially is 
$$
  0\ =\ \sum_\beta t^{\binom{n+1}{2}-|\beta|}
       L_\beta\, p_\beta(K)\;.
$$
For $K\in\Omega_\alpha$, the sum is over $\beta\leq\alpha$
by~(\ref{eq:sch-def}).
Removing the common factor 
$t^{\binom{n+1}{2}-|\alpha|}$ gives
 \begin{equation}\label{eq:red-eq}
  0\ =\ \sum_{\beta\leq\alpha} t^{|\alpha|-|\beta|}
       L_\beta\, p_\beta(K)\;.
 \end{equation}
The case $l_i=1$ of Theorem~\ref{thm:special} is implied by the case 
$\alpha=(n{-}k{+}1,\ldots, n{-}1, n)=\hat{1}$ 
($\Omega_{\hat{1}}=\Gr(k,n)$)
of the following theorem.

\begin{thm}[{\cite[Theorem 4.2]{So_trans}}]\label{thm:real-trans}
 Let $L\subset{\mathbb R}^n$ be a $(n-k)$-plane, none of whose Pl\"ucker
 coordinates vanishes.
 Then there exist real numbers $t_1,t_2,\ldots,t_{k(n-k)}\in{\mathbb R}^\times$
 such that for every $\alpha\in\binom{[n]}{k}$ the
 intersection 
$$
  \Omega_\alpha\;\cap\, \Omega(t_1.L)\,\cap\, \Omega(t_2.L)\,\cap\,
  \cdots\,\cap\,\Omega(t_{|\alpha|}.L)
$$ 
 is transverse (so it contains $d(\alpha)$ points) with all points real.
\end{thm}

\begin{proof}
We induct on $m$ to construct numbers
$t_1,t_2,\ldots,t_{k(n-k)}\in{\mathbb R}^\times$
having the property that, for all $\alpha\in\binom{[n]}{k}$ with $|\alpha|=m$, 
$$
  \Omega_\alpha\cap\Omega(t_1.L)\cap\Omega(t_2.L)\cap\cdots\cap\Omega(t_m.L)
$$
is transverse (over ${\mathbb C}$) and each of its $d(\alpha)$ points are
real.

The case $m=0$ is trivial, as $|\alpha|=0$ implies that 
$\alpha=\hat{0}$ and $\Omega_{\hat{0}}$ is a single
(real) point in $\Gr(k,n)$.
Suppose we have constructed $t_1,t_2,\ldots,t_m\in{\mathbb R}^\times$ 
with the above properties.

Let $\alpha\in\binom{[n]}{k}$ with $|\alpha|=m+1$ and consider
 \begin{equation}\label{eq:degen-sum}
  \left(\sum_{\beta\lessdot\alpha}\Omega_\beta\right) \cap\left(
    \Omega(t_1.L)\cap\cdots\cap\Omega(t_m.L)\right)\ =\ 
  \sum_{\beta\lessdot\alpha}\Omega_\beta
    \cap\Omega(t_1.L)\cap\cdots\cap\Omega(t_m.L)\,.
 \end{equation}
Each summand in the second sum is transverse (over ${\mathbb C}$) and
consists of $d(\beta)$ real points.
The intersection on the left will be transverse and consist of
$d(\alpha)=\sum_{\beta\lessdot\alpha}d(\beta)$ real points only if no two
summands in the second sum share a point in common.

If two summands $\beta,\gamma$ with $\beta\neq\gamma$ share a point, then 
 \begin{equation}\label{eq:empty}
  \Omega_\beta\,\cap\,\Omega_\gamma\,\cap\, 
  \Omega(t_1.L)\,\cap\,\Omega(t_2.L)\,\cap\,\cdots\,\cap\,\Omega(t_m.L)
 \end{equation}
is non-empty.
By~(\ref{eq:sch-def}), $\Omega_\beta\cap\Omega_\gamma$ is contained in a
union of Schubert varieties of dimension less than $m$.
By~(\ref{eq:red-eq}) the condition for a fixed
$k$-plane $K\in\Omega_\delta$ to meet $t.L$ non-trivially is a
polynomial of degree at most $|\delta|$ in $t$.
Thus a fixed $k$-plane $K\in\Omega_\delta$ lies in at most $|\delta|$ distinct 
$(n{-}k)$-planes in the family $t.L$.
Hence the intersection~(\ref{eq:empty}) is empty, and so the
summands in the second sum of~(\ref{eq:degen-sum}) are disjoint.
This argument also shows that the points in the summand indexed by $\beta$ lie
in $\Omega^\circ_\beta=\Omega_\beta-\bigcup_{\delta<\beta}\Omega_\delta$.

For $t\in{\mathbb R}$, let $Z_t\subset\Omega_\alpha$ be the set of $k$-planes
$K$ satisfying the polynomial~(\ref{eq:red-eq}).
For $t\neq 0$, we have $Z_t=\Omega_\alpha\cap\Omega(t.L)$.
Since $L$ has no vanishing Pl\"ucker coordinates, the constant term of that
polynomial is $p_\alpha(K)$, and so by Theorem~\ref{thm:Schubert-cover},
$$
  Z_0\ =\ \Omega_\alpha\cap H_\alpha\ =\ 
     \sum_{\beta\lessdot\alpha}\Omega_\beta\,.
$$
Thus $Z_0$ meets 
 \begin{equation}\label{eq:m-planes}
     \Omega(t_1.L)\,\cap\,\Omega(t_2.L)\,\cap\,\cdots\,\cap\,\Omega(t_m.L)
 \end{equation}
transversally (over ${\mathbb C}$) in $d(\alpha)$ real points.
We see that there exists $\epsilon_\alpha>0$  such that if
$0<t\leq\epsilon_\alpha$, 
then $Z_t$ meets~(\ref{eq:m-planes}) transversally (over ${\mathbb C}$) in
$d(\alpha)$ real points.
Since the intersection defining $Z_0$ is transverse along each
$\Omega^\circ_\beta$ for $\beta\lessdot\alpha$, 
we may assume that for $0<t\leq\epsilon_\alpha$,
$$
  \Omega_\alpha\,\cap\,\Omega(t.L)\,\cap\,
    \Omega(t_1.L)\,\cap\,\cdots\,\cap\,\Omega(t_m.L)
$$
is transverse (over ${\mathbb C}$) and consists of $d(\alpha)$ real points.

Let $t_{m+1}=\min\{\epsilon_\alpha\colon|\alpha|=m+1\}$.
Then for any $\alpha\in\binom{[n]}{k}$ with $|\alpha|= m+1$, 
$$
   \Omega_\alpha\,\cap\, \Omega(t_1.L)\cap\,
   \cdots\,\cap\,\Omega(t_m.L)\,\cap\,\Omega(t_{m+1}.L)
$$
is transverse (over ${\mathbb C}$) and consists of $d(\alpha)$ real points.
\end{proof}

\subsection{Further extensions of the Schubert calculus}

This special Schubert calculus admits further extensions, some of which are
known to be fully real.
Interestingly, some problems involving the Lagrangian
Grassmannian~\ref{sec:lagriangian} are known to be {\it fully unreal}, a new
phenomenon. 
Further investigation is encouraged; in particular, the Lagrangian Schubert
calculus may yield the first class of enumerative which are {\it not} fully
real.
\medskip

\subsubsection{General Schubert calculus}
A {\it flag} in $n$-space is a sequence of linear subspaces 
$\Fdot\colon F_1\subset F_2\subset\cdots\subset F_n$ with $\dim F_i=i$.
Given $\alpha\in\binom{[n]}{k}$, the Schubert condition of type $\alpha$ on a
$k$-plane $K$ imposed by the flag $\Fdot$ is
 \begin{equation}\label{eq:GenSchCond}
   \dim K\cap F_{n+1-\alpha_j}\ \geq\ k+1-j\qquad\mbox{for}\quad j=1,2,\ldots,n\,.
 \end{equation}
The Schubert variety $X_\alpha\Fdot\subset\Gr(k,n)$ is the set of all
$k$-planes $K$ satisfying~(\ref{eq:GenSchCond}).

We relate this to the definitions of Section~\ref{sec:SpSchCalc}.
Let $e_1,e_2,\ldots,e_n$ be a basis for ${\mathbb R}^n$ and for its
complexification ${\mathbb C}^n$.
Defining $F_i$ to be the linear span of $e_1,e_2,\ldots,e_i$ gives the 
flag $\Fdot$.
Then the Schubert variety $X_\alpha\Fdot$ is $\Omega_{\alpha^\vee}$, where
$\alpha^\vee_j=n+1-\alpha_{k+1-j}$ for each $j$, and so 
the codimension of $X_\alpha\Fdot$ is $|\alpha|$.
A special Schubert condition is when 
$\alpha=(1,\ldots,k{-}1,k{+}l)$
so that $X_\alpha\Fdot=\Omega(F_{n-k+1-l})$.
\smallskip

The general problem of the classical Schubert calculus of enumerative geometry
asks, given $\alpha^1,\alpha^2,\ldots,\alpha^s\in\binom{[n]}{k}$ with 
$|\alpha^1|+|\alpha^2|+\cdots+|\alpha^s|=k(n-k)$ and general flags
$\Fdot^1,\Fdot^2,\ldots,\Fdot^s\subset {\mathbb C}^n$, how many points are
there in the intersection\frankfootnote{In this survey
flags are general when the corresponding intersection is
transverse.} 
 \begin{equation}\label{eq:GenGrassInt}
   X_{\alpha^1}\Fdot^1\cap X_{\alpha^2}\Fdot^2\cap\cdots\cap
   X_{\alpha^s}\Fdot^s\,?
 \end{equation}
There are algorithms due to Pieri~\cite{Pi1893} and Giambelli~\cite{Gi1902} to
compute these numbers.
Other than the case when the $\alpha^i$ are indices of special Schubert
varieties, it remains open whether the general Schubert calculus
is fully real.
(See~\cite{So00b} and~\cite{So_shapiro-www} for some cases.)\frankfootnote{While
this survey was in review, Ravi Vakil communicated to the author a proof that
the classical Schubert calculus is fully real~\cite{Vakil}.} 
\medskip

\subsubsection{Quantum Schubert calculus}\label{sec:quantum}
The space $\calM^q_{k,n}$ of degree $q$ maps 
$M\colon{\mathbb P}^1\to\Gr(k,n)$ is a smooth quasi-projective
variety~\cite{Clark}.
A point $t\in{\mathbb P}^1$ and a Schubert variety $X_\alpha\Fdot$
together impose a quantum Schubert condition on maps 
$M\in\calM^q_{k,n}$,
$$
  M(t)\ \in\ X_\alpha\Fdot\,.
$$
The set of such maps has codimension $|\alpha|$ in $\calM^q_{k,n}$.
The quantum Schubert calculus of enumerative geometry asks the
following question.

\begin{ques}
 Given $\alpha^1,\alpha^2,\ldots,\alpha^s\in\binom{[n]}{k}$ with 
 $|\alpha^1|+|\alpha^2|+\cdots+|\alpha^s|=\dim\calM^q_{k,n}=qn+k(n-k)$, how many
 maps $M\in\calM^q_{k,n}$ satisfy
 $$
   M(t_i)\ \in\ X_{\alpha^i}\Fdot^i\qquad\mbox{for }i=1,2,\ldots,s\,,
 $$
 where $t_1,t_2,\ldots,t_s\in{\mathbb P}^1$ are general and
 $\Fdot^1,\Fdot^2,\ldots,\Fdot^s$ are general flags?
\end{ques}

Algorithms to compute this number were proposed by Vafa~\cite{Vafa} and
Intriligator~\cite{In91} and proven by 
Siebert and Tian~\cite{ST97} and by Bertram~\cite{Be97}.
A simple quantum Schubert condition defines a subvariety of codimension 1, 
 \begin{equation}\label{eq:simQSch}
   M(t)\ \cap\ L\ \neq\ \{0\}\,,
 \end{equation}
where $L$ is a $(n{-}k)$-plane.
Let $d(q;k,n)$ be the number of maps $M\in\calM^q_{k,n}$ satisfying
$\dim\calM^q_{k,n}$--many general simple quantum Schubert
conditions~(\ref{eq:simQSch}).
A combinatorial formula for this number was given by
Ravi, Rosenthal, and Wang~\cite{RRW98}.
For a survey of this particular enumerative problem and its importance to
linear systems theory, see~\cite{So01a}.

\begin{thm}[{\cite[Theorem 1.1]{So00a}}]
 Let $q\geq 0$ and $n>k>0$ and set $N:=\dim{\mathcal M}_{k,n}$.
 Then there exist $t_1,t_2,\ldots,t_N\in{\mathbb{RP}}^1$ and $(n{-}k)$-planes
 $L_1,L_2,\ldots,L_N\subset{\mathbb R}^n$ such that there are exactly
 $d(q;k,n)$ maps $M\in\calM^q_{k,n}({\mathbb C})$ satisfying
$$
  M(t_i)\cap L_i\ \neq\ \{0\}\qquad\mbox{for }i=1,2,\ldots,N\,,
$$
 and all of them are real.
\end{thm}

As with the classical Schubert calculus, the question of whether the general
quantum Schubert calculus is fully real remains open.
\medskip

\subsubsection{Schubert calculus of flags}\label{sec:flags}
Let ${\bf a}:=0<a_1<\cdots<a_r<a_{r+1}=n$ be a sequence of integers.
The manifold of partial flags in $n$-space (or the flag manifold) $\Fla$
is the collection of partial flags of subspaces 
$$
   \Edot\ :\ E_{a_1}\subset E_{a_2}\subset\cdots\subset E_{a_r}
$$
in $n$-space, where $\dim E_{a_i}=a_i$.
A complete flag $\Fdot$ is when ${\bf a}=1,2,\ldots,n{-}1$.

The Schubert varieties of $\Fla$ are indexed by permutations $w\in\calS_n$,
the symmetric group on $n$ letters, whose descents only occur at positions 
in ${\bf a}$.
That is, $w(i)>w(i+1)$ implies that $i\in\{a_1,a_2,\ldots,a_r\}$.
Let $I_{\bf a}$ be this set of indices.
For a complete flag $\Fdot$ and $w\in I_{\bf a}$, $\Fla$ has a Schubert variety
$$
  X_w\Fdot\ :=\ \left\{\Edot\in\Fla\mid \dim E_{a_j}\cap F_i\geq
  \#\{l\leq a_j\mid w(l)\leq i\}\right\}\,.
$$
The codimension of $X_w\Fdot$ is $\ell(w):=\#\{i<j\mid w(i)>w(j)\}$.
We state the general question in the Schubert calculus of enumerative geometry
for flags. 

\begin{ques}
 Given permutations $w_1,w_2,\ldots, w_s\in I_{\bf a}$ with
 $\ell(w_1)+\ell(w_2)+\cdots+\ell(w_s)=\dim\Fla$ and general flags
 $\Fdot^1,\Fdot^2,\ldots,\Fdot^s$, what is the number of points in the intersection
$$
    X_{w_1}\Fdot^1\cap X_{w_2}\Fdot^2 \cap\cdots\cap X_{w_s}\Fdot^s\ ?
$$
\end{ques}

There are algorithms~\cite{BGG73,De74} for computing this number and the 
numbers for the Lagrangian Schubert calculus in Section~\ref{sec:lagriangian}. 
When $w=(i,i+1)$, $X_w\Fdot$ is the simple Schubert variety, written
$X_i\Fdot$,
 \begin{equation}\label{eq:simpleSchVar}
   X_i\Fdot\ =\ X_{(i,i+1)}\Fdot\ =\ 
   \{\Edot\in\Fln\mid E_i\cap F_{n-i}\neq\{0\}\}\,.
 \end{equation}
\begin{thm}[{\cite[Corollary 2.2]{So00c}}]
 Given a list $i_1,i_2,\ldots,i_N$ ($N=\dim\Fla$) of numbers with 
 $i_j\in{\bf a}$, there exist real flags
 $\Fdot^1,\Fdot^2,\ldots,\Fdot^N$ such that the intersection of simple
 Schubert varieties
$$ 
  X_{i_1}\Fdot^1\cap X_{i_2}\Fdot^2\cap\cdots\cap  
  X_{i_N}\Fdot^N
$$
 is transverse and consists only of real flags.
\end{thm}

The case of ${\bf a}=2<n{-}2$ of this theorem was proven 
earlier~\cite[Theorem 13]{So97b}.
It remains open whether the general Schubert calculus of flags is fully real. 
 \medskip

\subsubsection{Unreality in the Lagrangian Schubert
calculus}\label{sec:lagriangian} 
Let $V$ be a $2n$-dimensional vector space equipped with the alternating
bilinear form
 \begin{equation}\label{eq:alt-form}
  \left\langle \sum x_ie_i\,,\ \sum y_je_j\right\rangle\ :=\ 
  \sum_{i=1}^n x_iy_{2n+1-i}-y_ix_{2n+1-i}\,.
 \end{equation}
A subspace $H\subset V$ is {\it isotropic} if the restriction of the form to
$H$ is identically zero, $\langle H,H\rangle\equiv 0$.
The dimension of an isotropic subspace is at most $n$ and {\it Lagrangian
subspaces} are isotropic subspaces with this maximal dimension.
The Lagrangian Grassmannian  $LG(n)$ is the set of Lagrangian subspaces of $V$,
an algebraic manifold of dimension $\binom{n+1}{2}$.

A flag $\Fdot$ in $V$ is {\it isotropic} if $F_n$ is Lagrangian and 
$\langle F_i,F_{2n-i}\rangle\equiv 0$ for all $i=1,2,\ldots,2n{-}1$.
Given an isotropic flag $\Fdot$, the Lagrangian Grassmannian has Schubert
varieties $\Psi_\lambda\Fdot$ indexed by decreasing sequences
$\lambda\colon n\geq \lambda_1>\lambda_2>\cdots>\lambda_l>0$ of positive integers,
called strict partitions. (Here $l$ can be any integer between 0 and $n$).
The codimension of $\Psi_\lambda\Fdot$ is
$|\lambda|=\lambda_1+\lambda_2+\cdots+\lambda_l$. 
The Lagrangian Schubert calculus asks for the number of points in a
transverse zero-dimensional intersection of Schubert varieties.

The {\it simple Schubert variety} $\Psi_1\Fdot$ consists of those Lagrangian
subspaces meeting the Lagrangian subspace $F_n$ non-trivially.
The simple Lagrangian Schubert calculus is fully unreal.

\begin{thm}\label{thm:no-real}
 There exist isotropic real flags
 $\Fdot^1,\Fdot^2,\ldots,\Fdot^{\binom{n+1}{2}}$ 
 such that the intersection of Schubert varieties 
$$
  \Psi_1\Fdot^1\cap\Psi_1\Fdot^2\cap\cdots\cap\Psi_1\Fdot^{\binom{n+1}{2}}
$$
 is transverse with {\bf no} real points.
\end{thm}

\begin{remark}
 In~\cite[Theorem 4.2]{So00c} flags were given so that the intersection had no
 real points, and it was not known if that intersection was transverse.
 Perturbing those flags slightly (so that the intersection becomes transverse)
 gives the above result.

 We do not know if these (or many other) enumerative problems in the
 Lagrangian Schubert calculus are fully real.
 Experimentation\frankfootnote{This will be reported in the forthcoming
 article~\cite{Sottile_Exp}.} suggests that the situation is complicated.
 Briefly, while many other enumerative problems in the Lagrangian Schubert
 calculus are fully unreal, there are a few which are fully real.
 For example, there 
 exist 2 real isotropic 2-planes and 2 real isotropic 3-planes such that
 all 4 Lagrangian subspaces meeting each of these are real (see
 Theorem~\ref{thm:Lagr-real-calc}). 

 The problems of Theorem~\ref{thm:no-real} may give examples
 of enumerative problems that are not fully real.
 Experimental evidence suggests however that this may be unlikely.
 For example, (case $n=3$ of Theorem~\ref{thm:no-real}) there are 16
 Lagrangian subspaces in ${\mathbb C}^6$ having 
 non-trivial intersection with 6 general Lagrangian subspaces.
 We computed 30,000 random instances of this enumerative problem, 
 and found two examples of 6 real Lagrangian subspaces such that all 16
 Lagrangian subspaces meeting them are real.
 Table~\ref{table:lagr-random} shows the number of these 30,000 systems having a
 given number of real solutions.
%
%
\begin{table}[htb]
 \begin{tabular}{|c||c|c|c|c|c|c|c|c|c|}\hline
  Number of real solutions&0&2&4&6&8&10&12&14&16\rule{0pt}{13pt}\\\hline
  Frequency&
   4983 & 8176 & 9314 & 5027 & 1978 & 445 & 67 & 8 & 2\rule{0pt}{13pt}\\\hline
 \end{tabular}\vspace{.1in}
 \caption{Frequency with given number of real
  solutions.\label{table:lagr-random}} 
\end{table}
\end{remark}

\section{The Conjecture of Shapiro and Shapiro}\label{sec:Shapiro}
 
The results of Section~\ref{sec:SchubertCalculus} were inspired by a
remarkable conjecture of Boris Shapiro and Michael Shapiro.
Let $\gamma\colon{\mathbb R}\to{\mathbb R}^n$ be the rational normal curve
$$
  \gamma(t)\ =\ (1,\, t,\, t^2,\,\ldots,\, t^{n-1})\ =\ 
                \sum_{i=1}^n t^{i-1}e_i\,.
$$
For $t\in{\mathbb R}$, define the flag $\Fdot(t)$ by any of the three
equivalent ways
 \begin{eqnarray}
  F_i(t)&=& \mbox{ linear span of } \gamma(t),\,
     \frac{d}{dt}\gamma(t),\,\ldots,\,\frac{d^{i-1}}{dt^{i-1}}\gamma(t)\,,
       \label{eq:5.1} \\
   &=& \mbox{ row space }\rule{0pt}{40pt}\left[ \begin{array}{ccccccc}
        1&t&t^2&&\cdots&&t^{n-1}\\
        0&1&2t&&\cdots&&(n{-}1)t^{n-2}\\
      \vdots&&\ddots&\ddots&&&\vdots\\
        0&\cdots&0&1&it&\cdots&\binom{n-1}{i-1}t^{n-i}
      \end{array}\right]\ , \label{eq:osc-matrix}\\
   &=& \mbox{ $i$-plane osculating rational normal curve $\gamma$ at
  $\gamma(t)$}\,.\rule{0pt}{20pt}
 \end{eqnarray}
This makes sense for $t\in{\mathbb C}$ and is extended to
$t=\infty\in{\mathbb P}^1$ by setting $F_i(\infty)$ to be the row space of the
last $i$ rows of the $n\times n$ identity matrix.

\begin{conj}[Shapiro-Shapiro]\label{conj:ShSh}
 Let $\alpha^1,\alpha^2,\ldots,\alpha^s\in\binom{[n]}{k}$ be such that
 $|\alpha^1|+|\alpha^2|+\ldots+|\alpha^s|=k(n-k)$.
 Then, for every distinct $t_1,t_2,\ldots,t_s\in{\mathbb R}$,
 the intersection of Schubert varieties
 \begin{equation}\label{eq:Shap-int}
  X_{\alpha^1}\Fdot(t_1)\cap X_{\alpha^2}\Fdot(t_2)\cap\cdots\cap
  X_{\alpha^s}\Fdot(t_s)
 \end{equation}
 is {\rm (a)} transverse, and\/ {\rm (b)} consists only of real points.
\end{conj}

Eisenbud's and Harris's dimensional transversality 
result~\cite[Theorem 2.3]{EH83} guarantees that the 
intersection~(\ref{eq:Shap-int}) is zero-dimensional. 
Not only does Conjecture~\ref{conj:ShSh} state that the classical Schubert
calculus is fully real, but it also proposes flags witnessing this full reality.
This conjecture has been central to subsequent developments in the real
Schubert calculus and it has direct connections to other parts of
mathematics, including linear systems theory and linear series on 
${\mathbb P}^1$ (see Remark~\ref{rem:MIC}). 
The article~\cite{So00b} and the web page~\cite{So_shapiro-www} give a more
complete discussion. 

One aspect of this conjecture which we relate is the following.

\begin{thm}[{\cite[Theorem 3.3]{So00b}}]\label{thm:implies}
 For a given $k$ and $n$, the general case of Conjecture~\ref{conj:ShSh}
 follows from the special case when each Schubert condition is simple, that is, 
 when each $\alpha^i=(1,2,\ldots,k{-}1,k{+}1)$.
\end{thm}

Consider these osculating flags $\Fdot(t)$ in more detail.
If the $i$th row of the matrix~(\ref{eq:osc-matrix})
is multiplied by $t^i$ (which does not affect its row space when $t\neq 0$),
then the entry in position $i,j$ is $\binom{j-1}{i-1}t^j$, and so we have
$$
   F_i(t)\ =\ t.F_i(0)\,,
$$
where $t.F_i(0)$ is given by the action~(\ref{eq:action}) of ${\mathbb R}^\times$
on ${\mathbb R}^n$.
The $\alpha$th Pl\"ucker coordinate of $F_i(0)$ is
\begin{equation}\label{eq:pl-form}
  p_\alpha( F_i(0))\ =\ \prod_{j<l}\frac{\alpha_j-\alpha_l}{j-l}\,,
\end{equation}
which is non-vanishing.
Thus Theorem~\ref{thm:real-trans} has the following corollary.

\begin{thm}[{\cite[Theorem 1]{So99a}}]\label{thm:osc-flags}
 There exist $t_1,t_2,\ldots,t_{k(n-k)}\in{\mathbb R}$ such that 
 there are exactly $d(k,n)$ $k$-planes meeting each $(n{-}k)$-plane
 $F_{n-k}(t_i)$ non-trivially, and all are real.
 Equivalently, if $\alpha=(1,2,\ldots,k{-}1,k{+}1)$ so that $|\alpha|=1$, then
 the intersection of Schubert varieties
 \begin{equation}\label{eq:shap-int}
   X_\alpha\Fdot(t_1)\cap X_\alpha\Fdot(t_2)\cap \cdots\cap 
   X_\alpha\Fdot(t_{k(n-k)})
 \end{equation}
 is transverse with all points real.
\end{thm}

This establishes a weak form of Conjecture~\ref{conj:ShSh} for simple Schubert
conditions, replacing the quantifier {\it for all} $t_i\in{\mathbb R}$ by
{\it there exists} $t_i\in{\mathbb R}$.

If the parameters $t_i$ in~(\ref{eq:shap-int}) vary, then the number of real
points in that intersection could change, but only if two points first collide
(prior to spawning a complex conjugate pair of solutions).
This is the reverse of the progression in Dietmaier's algorithm, as displayed in
Figure~\ref{fig:Dietmaier}. 
This situation cannot occur if the intersection remains transverse.
Together with Theorems~\ref{thm:implies} and~\ref{thm:osc-flags}
we deduce the following result.

\begin{thm}[{\cite[Theorem 6]{So99a}}]
 Part {\rm (a)} of Conjecture~\ref{conj:ShSh} for simple conditions implies 
 part {\rm (b)} for arbitrary Schubert conditions.
\end{thm}

If $t>0$, then by~(\ref{eq:pl-form}) the Pl\"ucker coordinates of $F_i(t)$ are
strictly positive.
An upper triangular $n\times n$ matrix $g$ is {\it totally positive} if every
$i\times i$ subdeterminant of $g$ is non-negative, and vanishes only if that
subdeterminant vanishes on all upper triangular matrices.
For example, when $t>0$ and $i=n$, the matrix~(\ref{eq:osc-matrix}) is totally
positive.
Write $\Gamma(t)$ for this matrix.
It has the form $e^{t\eta}$ where $\eta$ is the principal nilpotent matrix 
$\Gamma'(0)$.
Observe that if $t_1<t_2<\cdots<t_s$, then 
$$
  \Fdot(t_i)\ =\ \Gamma(t_i-t_{i-1}) \cdot \Fdot(t_{i-1})
$$
Conjecture~\ref{conj:ShSh} has a more general version involving totally
positive matrices.

\begin{conj}[Shapiro-Shapiro~{\cite[Conjecture 4.1]{So00b}}]
 Let $\alpha^1,\alpha^2,\ldots,\alpha^s\in\binom{[n]}{k}$ be such that
 $|\alpha^1|+|\alpha^2|+\ldots+|\alpha^s|=k(n-k)$ and 
 suppose $g_2,g_3,\ldots,g_s$ are totally positive matrices.
 Given any real flag $\Fdot^1$, define $\Fdot^i$ for $i>1$ by
 $\Fdot^i:=g_i\cdot\Fdot^{i-1}$,
 the intersection of Schubert varieties
$$
  X_{\alpha^1}\Fdot^1\,\cap\, X_{\alpha^2}\Fdot^2\,\cap\,\cdots\,\cap\,
  X_{\alpha^s}\Fdot^s
$$
 is transverse with all points real.
\end{conj}

There is some experimental evidence for this version of
Conjecture~\ref{conj:ShSh}. 
Subsequent conjectures involving osculating flags have
versions involving totally positive matrices.
We leave their statements to the reader, they will be explored more fully
in~\cite{Sottile_Exp}. 

\subsection{Rational functions with real critical points.}\label{eq:ratreal}
By far the strongest evidence for
Conjecture~\ref{conj:ShSh}
is that it is true when $k$ or $n-k$ is equal to 2.

\begin{thm}[Eremenko and Gabrielov~\cite{EG00}]
 Conjecture~\ref{conj:ShSh} is true when one of $k$ or $n{-}k$ is $2$.
\end{thm}

This is a consequence of a theorem about rational functions with real
critical points.
A rational function is an algebraic map
$\varphi\colon{\mathbb P}^1\to{\mathbb P}^1$.
Two rational functions $\varphi_1$ and $\varphi_2$ are equivalent if
$\varphi_1=g\circ\varphi_2$, where $g$ is a fractional linear transformation.

\begin{thm}[\cite{EG00}]\label{thm:ratfn}
 If all the critical points of a rational function are real, then it is
 equivalent to a real rational function.
\end{thm}

Consider the composition
 \begin{equation}\label{eq:ratmap}
  {\mathbb P}^1\ \longrightarrow\ 
  {\mathbb P}^d\ \relbar\to\ {\mathbb P}^1\,,
 \end{equation}
where the first map is the rational normal curve
$$
  \gamma\ \colon\ [s,t]\ \longmapsto\ 
    [s^d,\,s^{d-1}t,\,\ldots,\,st^{d-1},\,t^d]\,,
$$
and the second is a linear projection
$$
  [x_0,x_1,\ldots,x_d]\ \longmapsto\ [\Lambda_1(x),\Lambda_2(x)]\,,
$$
where $\Lambda_1$ and $\Lambda_2$ are independent linear forms.
Let $E\subset{\mathbb P}^d$ be the center of this projection,
the linear subspace where $\Lambda_1=\Lambda_2=0$.
When $E$ is disjoint from the rational normal curve, this composition defines
a rational function of degree $d$, and all such rational functions occur in
this manner. 
In fact, equivalence classes of rational maps are exactly those maps with the
same center of projection.
A rational function $\varphi$ has a critical point at $t\in{\mathbb P}^1$ if
$d\varphi$ vanishes at $t$.
If we consider the composition~(\ref{eq:ratmap}), this implies that the center
$E$ meets the line tangent to the rational normal curve 
$\gamma({\mathbb P}^1)$ at $\gamma(t)$.

Goldberg~\cite{Go91} asked (and answered) the question:
how many equivalence classes of rational functions of degree $d$ have a given
set of $2d-2$ critical points?
Reasoning as above, she reduced this to the problem of determining the number
of codimension 2 planes $E$ in ${\mathbb P}^d$ meet $2d-2$ given tangents to
the rational normal curve.
Formulating this in the dual projective space, we recover the problem
of the Introduction to Section~\ref{sec:SchubertCalculus}:
Determine the number of lines in ${\mathbb P}^d$ meeting $2d-2$ general
codimension 2 planes.
The answer is the $d$th Catalan number, 
$C_d=\frac{1}{d}\binom{2d-1}{d-1}$~(\ref{eq:Catalan}). 

Consider the above description
in the affine cone over ${\mathbb P}^d$.
The centers $E\in\Gr(d{-}1,d{+}1)$ giving a rational function of degree $d$
with critical point at $t\in{\mathbb C}\subset{\mathbb P}^1$ are points 
in the Schubert variety $X_{\alpha}\Fdot(t)$, where $\alpha$ is 
the simple Schubert condition of Theorem~\ref{thm:implies} and $\Fdot(t)$
is the flag of subspaces osculating the rational normal curve 
$\gamma({\mathbb P}^1)$ at $\gamma(t)$. 
An equivalence class of rational functions contains a real rational
function when the common center $E$ is real.
In this way, Theorem~\ref{thm:ratfn} implies
Conjecture~\ref{conj:ShSh} when $n{-}k=2$ and each $|\alpha^i|=1$.
By Theorem~\ref{thm:implies}, this implies the full conjecture when
$n{-}k=2$. 
Working in the dual space, we deduce Conjecture~\ref{conj:ShSh} when $k=2$.

Gabrielov and Eremenko prove Theorem~\ref{thm:ratfn} by showing there exist
$C_d$ distinct real rational functions with critical points at a given set
of $2d{-}2$ real numbers.
Let $\varphi\colon{\mathbb P}^1\to{\mathbb P}^1$ be a real rational function of
degree $d$ with only real critical points.
Then $\varphi^{-1}({\mathbb P}^1_{\mathbb R})\subset{\mathbb P}^1_{\mathbb C}$
is an embedded graph containing  
${\mathbb P}^1_{\mathbb R}$ which is stable under complex conjugation and
whose vertices have valence 4 and are at the critical points.
There are exactly $C_d$ isotopy classes of such graphs.
For each isotopy class $\Gamma$ and collection of $2d{-}2$ distinct real points,
they essentially construct a rational function $\varphi$ of degree $d$
having these critical points with 
$\varphi^{-1}({\mathbb P}^1_{\mathbb R}) \in\Gamma$.
A key point involves degenerate rational maps with
fewer than $2d{-}2$ critical points.
It may be interesting to relate this to the degenerations in the
proof of Theorem~\ref{thm:real-trans}.

\begin{remark}\label{rem:MIC}
 If we project the rational normal curve from a center $E$ of codimension $m+1$,
 then the image of the projection is ${\mathbb P}^m$ and the composition is a
 parameterized rational curve in ${\mathbb P}^m$.
 Post-composition by an element of Aut(${\mathbb P}^m$) defines an equivalence
 relation on such maps and equivalence classes are determined by the centers
 of projection.
 These centers are codimension $m+1$ linear series of degree $d$ on 
 ${\mathbb P}^1$.
 Such a linear series is {\it ramified} at a point $t\in{\mathbb P}^1$ when
 the center $E$ meets the $m$-plane osculating the rational normal curve.
 A rational curve/linear series is {\it maximally inflected}\/ if  the
 ramification is at real points.
 We just considered the case $d=1$ of rational functions with real critical
 points, and the case $d=2$ was discussed in Section~\ref{old-real}.
 The existence of maximally inflected curves with simple ramification is a
 consequence of Theorem~\ref{thm:osc-flags}, and Conjecture~\ref{conj:ShSh}
 predicts this is a rich class of real curves in ${\mathbb P}^m$.
 This connection between linear series and the Schubert calculus originated in 
 work of Castelnuovo~\cite{Ca1889}.
\end{remark}

\subsection{Generalizations of Conjecture~\ref{conj:ShSh}}
The Grassmannian, flag manifolds, and Lagrangian Grassmannian 
are examples of flag varieties $G/P$ where $G$ is a reductive algebraic group
and $P$ a parabolic subgroup.
These flag varieties have Schubert varieties and the most general form of the
Schubert calculus involves zero-dimensional intersections of these Schubert
varieties. 
Likewise, these flag varieties have real forms (given by
split\frankfootnote{Split is a technical term: 
${\mathbb R}^\times\subset{\mathbb C}^\times$ is a split form of $GL_1$, but 
$S^1\subset{\mathbb C}^\times$ is not.}
forms of $G$ and $P$) and there is a generalization of
Conjecture~\ref{conj:ShSh} for these real forms.
This generalization is false, but in a very interesting way.
We describe what is known about this general conjecture for the flag
manifolds and the Lagrangian Grassmannian, and give conjectures
describing what we believe to be true.
\medskip

\subsubsection{The manifold of partial flags}
Let ${\bf a}:=0<a_1<a_2<\cdots<a_r<a_{r+1}=n$ be a sequence of integers.
The straightforward generalization (which was its original form and which is
false) of Conjecture~\ref{conj:ShSh} for the flag manifold $\Fla$ declares
that a zero-dimensional 
intersection of Schubert varieties given by flags osculating the real rational
normal curve consists only of real points.

For $\alpha\in\binom{[n]}{a_i}$, we have the {\it Grassmannian Schubert variety}
$$
  Y_{\alpha,i}\Fdot\ :=\ \{\Edot\in\Fla\mid
    E_{a_i}\in X_\alpha\Fdot\subset \Gr(a_i,n)\}\,.
$$
For example, if $|\alpha|=1$, then 
$Y_{\alpha,i}\Fdot$ is the simple Schubert variety $X_i\Fdot$ of
Section~\ref{sec:flags}.

\begin{example}
 Let ${\bf a}=2<3<5=n$ so that $\Fla$ is the manifold of partial flags 
 $E_2\subset E_3$ in 5-space.
 This flag manifold has two simple Schubert varieties
 \begin{eqnarray*}
  X_2\Fdot &=& \{E_2\subset E_3\mid E_2\cap F_3\neq\{0\}\}\\
  X_3\Fdot &=& \{E_2\subset E_3\mid E_3\cap F_2\neq\{0\}\}
 \end{eqnarray*}
A calculation~\cite[Example 2.5]{So00c} shows that
\begin{equation}\label{eq:no-real}
  \begin{array}{c}
\quad  X_2\Fdot (-8)\,\cap\, X_3\Fdot(-4)\,\cap\, X_2\Fdot(-2)\,\cap 
            X_3\Fdot(-1)\,\cap\,\vspace{4pt}\\
   \,X_2\Fdot(1)\,\cap\,  X_3\Fdot(2)\,\cap\, X_2\Fdot(4)\cap\,
              X_3\Fdot(8)
 \end{array}
\end{equation}
is transverse with {\it none} of its 12 points real.

Thus the straightforward generalization of Conjecture~\ref{conj:ShSh} is
completely false.
On the other hand, if $t_1<t_2<\cdots<t_8$ are any of the 24,310 subsets of
eight numbers from  
$$
  \{-6,-5,-4,-3,-2,-1,1,2,3,5,7,11,13,17,19,23,29\}\,,
$$
then the intersection 
\begin{equation}\label{eq:gen12int}
  \begin{array}{c}
  \quad X_2\Fdot (t_1)\,\cap\, X_2\Fdot(t_2)\,\cap\, 
              X_2\Fdot(t_3)\,\cap\, X_2\Fdot(t_4)\vspace{4pt}\\
  \cap\, X_3\Fdot(t_5)\,\cap\, X_3\Fdot(t_6)\,\cap\,
                 X_3\Fdot(t_7\,)\cap\, X_3\Fdot(t_8)
 \end{array}
\end{equation}
is transverse with {\it all} of its 12 points real~\cite[Example 2.5]{So00c}.
\end{example}

We have used more than $2\times 10^7$ seconds of CPU time 
%
%
investigating this problem of intersections of Schubert varieties in manifolds
of partial flags given by flags osculating the rational normal curve, and a
picture is emerging of what to expect, at least for Grassmannian Schubert
varieties.

Suppose we have a list of indices of Grassmannian Schubert varieties
$$
  (\alpha^1,i_1),\,(\alpha^2,i_2),\,\ldots,\,(\alpha^s,i_s),\qquad
   \mbox{with}\quad \alpha^j\in\binom{[n]}{i_j}
$$
where $i_j\in\{a_1,\ldots,a_r\}$ and 
$|\alpha^1|+|\alpha^2|+\cdots+|\alpha^s|=\dim\Fla$.
Call such a list {\it Grassmannian Schubert data} for $\Fla$.
Consider an intersection of Grassmannian Schubert varieties
 \begin{equation}\label{eq:GSVInt}
  Y_{\alpha^1,i_1}\Fdot(t_1)\cap   Y_{\alpha^2,i_2}\Fdot(t_2)\cap 
   \cdots\cap Y_{\alpha^s,i_s}\Fdot(t_s)\,,
 \end{equation}
where $t_1 <t_2<\cdots<t_s$ are distinct real numbers and
$\Fdot(t)$ is the flag osculating the rational normal curve $\gamma$ at
$\gamma(t)$.

\begin{conj}\label{conj:flags}
 Let\/ ${\bf a}=1< a_1<a_2<\cdots<a_r<a_{r+1}=n$ and 
$$
  (\alpha^1,i_1),\ (\alpha^2,i_2),\ \ldots,\ (\alpha^s,i_s),\qquad
  \mbox{with}\quad \alpha^j\in\binom{[n]}{i_j}
$$
 be Grassmannian Schubert data for\/ $\Fla$. 
\begin{enumerate}
  \item
       For every choice of real numbers
       $t_1<t_2<\cdots< t_s$, the intersection~$(\ref{eq:GSVInt})$ 
       is {\rm (a)} transverse with 
       {\rm (b)} with all points of intersection real if the indices $i_j$ are 
       in order:
       $i_1\leq i_2\leq\cdots\leq i_s$ or else
       $i_1\geq i_2\geq\cdots\geq i_s$.
  \item
       If the indices $i_j$ are not in order, nor is any cyclic permutation of
       the indices, 
       then there exist real numbers $t_1<t_2<\cdots< t_s$ such that the
       intersection~$(\ref{eq:GSVInt})$ is transverse with {\bf all} points
       real, and there exist real numbers $t_1<t_2<\cdots< t_s$ 
       such that the intersection is transverse with {\bf not all} points real.
\end{enumerate}
 In particular, enumerative problems involving Grassmannian Schubert varieties
 on $\Fla$ are fully real.
\end{conj}

\begin{remark}\mbox{ }
 \begin{enumerate}
  \item
       When each $|\alpha^i|=1$, there are real numbers 
       $t_1<t_2<\ldots<t_s$ such that the intersection~(\ref{eq:GSVInt}) is
       transverse with all points real~\cite[Corollary 2.2]{So00c}.
       Thus the number of real points in the intersection~(\ref{eq:GSVInt})
       is expected to vary when the the indices $i_j$ are not in order.
  \item
       When $r=1$, $\Fla$ is the Grassmannian $\Gr(a_1,n)$ and the 
       condition on the indices being ordered is empty,
       so this case of Conjecture~\ref{conj:flags} reduces to
       Conjecture~\ref{conj:ShSh}. 
  \item
       There is considerable evidence for this conjecture when $r=2$.
       Part (2) is true for every set of Grassmannian Schubert
       data in ${\mathbb F}\ell_{2<3}{\mathbb C}^5$ and 
       all except one such set in 
       ${\mathbb F}\ell_{2<4}{\mathbb C}^6$\frankfootnote{The 
       one exception is due to our inability to compute any instances.}.
       Many instances of these same enumerative problems with ordered
       partitions have been computed, and in each instance of (1) the
       intersection~(\ref{eq:GSVInt}) is transverse with all points real.
  \item
       As for Conjecture~\ref{conj:ShSh}, (a) implies (b) in Part
       (1) of Conjecture~\ref{conj:flags}.
  \item
       We have tested no instances of the intersection~(\ref{eq:GSVInt}) with
       $r>2$, so the truth may differ from the exact statement of
       Conjecture~\ref{conj:flags}.  
  \item
       Conjecture~\ref{conj:flags} has nothing to say when the Schubert data
       are not Grassmannian.
 \end{enumerate}
\end{remark}

\medskip

\subsubsection{The Lagrangian Grassmannian}
As in Section~\ref{sec:lagriangian}, let $V$ be a $2n$-dimensional vector
space equipped with the alternating form $\langle\cdot,\cdot\rangle$ defined
by~(\ref{eq:alt-form}). 
Set
 \begin{equation}\label{eq:isornc}
  \gamma(t)\ :=\ \left(
      1,\,t,\,\frac{t^2}{2},\,\ldots,\, \frac{t^n}{n!},\,
  -\frac{t^{n+1}}{(n+1)!},\,\frac{t^{n+2}}{(n+2)!},\,
      \ldots,\,(-1)^{n-1}\frac{t^{2n-1}}{(2n-1)!}\right)\,.
 \end{equation}
The flag $\Fdot(t)$ of subspaces osculating this
rational normal curve is isotropic for all $t\in{\mathbb C}$.

Given a strict partition $\lambda\colon\lambda_1>\lambda_2>\cdots>\lambda_l>0$
and an isotropic flag $\Fdot$, set
$$
  \Psi_\lambda\Fdot\ :=\ \{H\in LG(n)\mid H\cap F_{n+1-\lambda_i}\geq i
    \mbox{ for } i=1,\ldots,l\}\,.
$$
The codimension of this Schubert variety is
$|\lambda|=\lambda_1+\lambda_2+\cdots+\lambda_l$. 
For example, 
$$
  \Psi_a\Fdot=\{H\mid H\cap F_{n+1-a}\neq\{0\}\}\,.
$$
With these definitions, we state the version of Theorem~\ref{thm:no-real}
proven in~\cite{So00c}.

\begin{thm}[{\cite[Theorem 4.2]{So00c}}]
 Given a strict partition $\lambda$, let
 $N:=\binom{n+1}{2}-|\lambda|=\dim\Psi_\lambda\Fdot$.
 If $N>1$, then there exist distinct real numbers $t_1,t_2,\ldots,t_N$ such
 that
$$
  \Psi_\lambda\Fdot(0)\,\cap\,\Psi_1\Fdot(t_1)\,\cap\,\Psi_1\Fdot(t_2)
  \,\cap\,\cdots\,\cap\,\Psi_1\Fdot(t_N)
$$
 is zero-dimensional with {\bf no} points real.
\end{thm}

Thus the generalization of Conjecture~\ref{conj:ShSh} is badly false for the
Lagrangian Grassmannian.
On the other hand it holds for some enumerative problems in the Lagrangian
Grassmannian.

\begin{thm}[\cite{Sottile_Exp}]\label{thm:Lagr-real-calc}
 For any distinct real numbers $t_1,t_2,t_3,t_4$, the intersection of
 Schubert varieties in $LG(3)$
$$
  \Psi_1\Fdot(t_1)\,\cap\,\Psi_1\Fdot(t_2)\,\cap\,
  \Psi_2\Fdot(t_3)\,\cap\,\Psi_2\Fdot(t_4)
$$
 is transverse with {\bf all} points real.
\end{thm}

These two results illustrate a dichotomy that is emerging from
experimentation.
Call a list of strict partitions $\lambda^1,\lambda^2,\ldots,\lambda^s$ with
$|\lambda^1|+|\lambda^2|+\cdots+|\lambda^s|=\binom{n+1}{2}$ 
Lagrangian Schubert data.
In every instance we have computed of a zero-dimensional intersection of
Lagrangian Schubert varieties whose flags osculate the rational normal
curve~(\ref{eq:isornc}), the intersection has been transverse with either all
points real or no points real.
Most interestingly, the outcome---all real or no real---has depended only upon
the Lagrangian Schubert data of the intersection.

\begin{conj}
 Given Lagrangian Schubert data $\lambda^1,\lambda^2,\ldots,\lambda^s$ and
 distinct real numbers $t_1,t_2,\ldots,t_s$, the intersection
$$
  \Psi_{\lambda^1}\Fdot(t_1)\,\cap\,\Psi_2{\lambda^2}\Fdot(t_2)\,\cap\,
  \cdots\,\cap\,\Psi_{\lambda^s}\Fdot(t_s)
$$
 is transverse with either
\begin{enumerate}
 \item[(a)] all points real, or
 \item[(b)] no points real, 
\end{enumerate}
 and the outcome {\rm (a)} or {\rm (b)} depends only upon the list
 $\lambda^1,\lambda^2,\ldots,\lambda^s$. 
\end{conj}

We do not have a good idea what distinguishes the Lagrangian Schubert data
giving all points real from the data giving no points real.
Further experimentation is needed.
\medskip

\subsubsection{Further generalizations of Conjecture~\ref{conj:ShSh}}
 The status of the generalizations of Conjecture~\ref{conj:ShSh} to
 other flag varieties is almost completely unknown.
 There is one flag variety, the Orthogonal Grassmannian, for which much
 is known.
 In particular, the analog of Theorem~\ref{thm:osc-flags} (involving
 codimension 1 Schubert varieties given by isotropic flags osculating the
 rational normal curve) holds for the Orthogonal
 Grassmannian~\cite[Corollary 3.3]{So00c}.  
 There has also been a significant amount of computer experimentation testing
 cases of the obvious generalization of Conjecture~\ref{conj:ShSh} for the
 orthogonal Grassmannian, and in
 each, all points of intersection were found to be real.
 Lastly, we remark that we do not know of any reasonable version of
 Conjecture~\ref{conj:ShSh} for the quantum
 Schubert calculus.

\section{Lower Bounds in the Schubert calculus}\label{sec:lower}

In Section~\ref{sec:ratcubics}, we saw that of the 12 rational cubics meeting
8 real points in the plane, at least 8 were real.
This was the first instance of a non-trivial lower bound on the number of real
solutions to a problem in enumerative geometry.
Recent work of Eremenko and Gabrielov shows this phenomenon is pervasive in
the Schubert calculus.

Recall from Section~\ref{sec:SpSchCalc} that the $k$-planes in ${\mathbb C}^n$ 
meeting $k(n{-}k)$ general $(n{-}k)$-planes non-trivially is a complementary 
dimensional linear section of the Grassmannian,  
$$
  \Lambda\,\cap\,\Gr(k,n)\,,
$$
where $\Lambda$ has codimension $k(n{-}k)$ in Pl\"ucker space.  
The number of such $k$-planes is the degree of $\Gr(k,n)$.
This is also the degree of the linear projection
$$
  \pi_E\ :\ \Gr(k,n)\ \hookrightarrow\ {\mathbb P}^{\binom{n}{k}-1}\ 
    \relbar\to\ {\mathbb P}^{k(n-k)}
$$
with center of projection a plane $E$ with codimension $k(n{-}k)+1$ disjoint 
from $\Gr(k,n)$. 
The connection between these two definitions of degree is that when  
$E\subset\Lambda$,  $\pi_E(\Lambda)$ is a point  
$x\in{\mathbb P}^{k(n-k)}$ and    
$$
  \Lambda\,\cap\,\Gr(k,n)\ =\ \pi_E^{-1}(x)\,.
$$
Since complex manifolds are canonically oriented, the degree of such a
linear projection is just the number of points in the inverse image of a
regular value $x$.

An important such linear projection is the Wronski map.
Let $L(t)$ be the $(n{-}k)$-plane osculating the rational normal
curve~(\ref{eq:5.1}).
By~(\ref{eq:laplace}) and the discussion following Theorem~\ref{thm:implies},
the equation for a $k$-plane $K$ to meet $L(t)$ non-trivially is
 \begin{equation}\label{eq:wronski}
  \sum t^{\binom{n+1}{2}-\binom{n-k+1}{2}-|\alpha|}L_\alpha p_\alpha(K)\,,
 \end{equation}
where $p_\alpha(K)$ is the $\alpha$th Pl\"ucker coordinate of $K$ and
$L_\alpha$ is the appropriately signed Pl\"ucker coordinate of $L(0)$
complementary to $\alpha$.
The association of a $k$-plane $K$ to the polynomial~(\ref{eq:wronski}) is the
{\it Wronski map}
$$
  \pi_W\ :\ \Gr(k,n)\ \longrightarrow\ {\mathbb P}^{k(n-k)}\,,
$$
where ${\mathbb P}^{k(n-k)}$ is the space of polynomials of degree at most
$k(n{-}k)$, modulo scalars. 
This is a linear projection as the coefficients in~(\ref{eq:wronski}) are
linear in the Pl\"ucker coordinates.
If $f$ has distinct roots $t_1,t_2,\ldots,t_{k(n-k)}$, then 
$\pi^{-1}_W(f)$ is the set of $k$-planes meeting
each of $L(t_1),L(t_2),\ldots,L(t_{k(n-k)})$ non-trivially.
Observe that $\pi^{-1}_W(f)$ is given by real linear equations on the
Grassmannian if and only if $f$ has real coefficients, which includes the case
when the $t_i$ are all real (the situation of the Shapiro
Conjecture~\ref{conj:ShSh}).

To see why this is called the Wronski map, consider $\Gr(k,n)$ as the set
of $k$-planes in the space of polynomials of degree at most $n{-}1$.
Given a $k$-plane
$$
  K\ =\ \mbox{linear span}\{f_1,f_2,\ldots,f_k\}\,,
$$
the Wronski determinant of $K$ is
$$
  W(K)\ =\ \det\ \left[\begin{array}{cccc}
               f_1&f_2&\cdots&f_k\\
               f'_1&f'_2&\cdots&f'_k\\
              \vdots&\vdots&\ddots&\vdots\\
         f^{(k-1)}_1&f^{(k-1)}_2&\cdots&f^{(k-1)}_k
         \end{array}\right]\ ,
$$
a polynomial of degree at most $k(n{-}k)$, well-defined modulo scalars.
Under a choice of coordinates given by the coefficients of a polynomial,
$W(K)=\pi_W(K)$.\medskip

Since real manifolds are not necessarily orientable, the degree of a map is a
${\mathbb Z}/2{\mathbb Z}$-valued invariant.
However, Kronecker~\cite{Kr1968} defined the degree of a regular map
${\mathbb P}^2_{\mathbb R}\to{\mathbb P}^2_{\mathbb R}$, which he called the
{\it characteristic}, and his definition makes sense for many maps 
$f\colon X\to Y$ of (not necessarily orientable) compact manifolds.

First suppose that $X$ is oriented.
For a regular value $y\in Y$ of $f$, define
$$
  \Char\;f\ :=\ \left| \sum_{x\in f^{-1}(y)}
            \sgn\det df_x\,\right|\ ,
$$
using local coordinates in $X$ consistent with its
orientation and any local coordinate near $y\in Y$.
The sum is well-defined up to multiplication by $\pm 1$.
This number is independent of choices, if $Y$ is connected.

When $X$ may not be orientable, let $\widetilde{X}$ be the space of orientations
of $X$, which is canonically oriented, and similarly let $\widetilde{Y}$ be the
space of orientations of $Y$.
Then $\widetilde{X}\to X$  and $\widetilde{Y}\to Y$ are 2 to 1 coverings of $X$
and $Y$, with covering group ${\mathbb Z}/2{\mathbb Z}$.
The map $f$ is {\it orientable} if it has a lift
$\widetilde{f}\colon\widetilde{X}\to\widetilde{Y}$ that is equivariant with
respect to the covering group ${\mathbb Z}/2{\mathbb Z}$.
Define the characteristic of an orientable map $f$ to be
$\Char\widetilde{f}$.
This is well-defined, even if $\widetilde{Y}$ consists of 2 components.
This characteristic satisfies a fundamental inequality.

\begin{prop}
 Let $f\colon X\to Y$ be an orientable map.
 Then, for every regular value $y\in Y$,
$$
  \# f^{-1}(y)\ \geq\ \Char f\,.
$$
\end{prop}

We call this characteristic the {\it real degree} of the map $f$.

Consider this notion for the Wronski map of the real Grassmannian 
$\Gr_{\mathbb R}(k,n)$.
The Grassmannian has dimension $k(n{-}k)$ and is orientable if and only if $n$
is even. 
Eremenko and Gabrielov~\cite{EG02} compute the degree of the Wronski map.
To state their result, we introduce some additional combinatorics.
Recall the interpretation~(\ref{eq:deg-BO})
$$
  \mbox{degree }\Gr(k,n)\ =\ \#\mbox{ chains in Bruhat order from 
    $\hat{0}$ to $\hat{1}$}\,,
$$
where $\hat{1}=(n{-}k{+}1,\ldots,n{-}1,n)$ is the top element in the Bruhat
order. 
We introduce a statistic on these chains.
Each cover $\alpha\lessdot\beta$ in the Bruhat order has a unique index
$i$ with 
$$
  \alpha_i\ <\ \beta_i\ =\ \alpha_i+1\qquad
  \mbox{ but }\qquad \alpha_j\ =\ \beta_j\quad\mbox{for }j\ \neq\ i\,.
$$
The word $w(\CH)$ of a chain 
$\CH\colon\hat{0}\lessdot\alpha^1\lessdot\alpha^2\lessdot
                       \cdots\lessdot\alpha^{k(n-k)}=\hat{1}$
is the sequence of indices $(i_1,i_2,\ldots,i_{k(n-k)})$ of these indices for
covers in $\CH$.
An inversion in such a word is a pair $j<l$ with $i_j<i_l$, and the weight,
$\omega(\CH)$ is the number of inversions in the word of the chain $\CH$.
For instance the chain in the Bruhat order of $\Gr(3,6)$ highlighted in
Figure~\ref{fig:BrDeg} 
$$
  123\,\lessdot\, 124\,\lessdot\, 125\,\lessdot\, 135\,\lessdot\, 
  145\,\lessdot\, 245\,\lessdot\, 246\,\lessdot\, 346\,\lessdot\, 
  356\,\lessdot\, 456\,,
$$
has word $332213121$ and length 5.
Define the inversion polynomial
$$
  I(k,n)(q)\ :=\ \sum q^{\omega(\ch)}\,,
$$
the sum over all chains $\CH$ in the Bruhat order from $\hat{0}$ to $\hat{1}$.

\begin{thm}[Eremenko and Gabrielov~{\cite[Theorem 1]{EG02}}]
 The characteristic of the Wronski map 
 $\pi_W\colon \Gr_{\mathbb R}(k,n)\to{\mathbb P}_{\mathbb R}^{k(n-k)}$ is
 $|I(k,n)(-1)|$.
\end{thm}

They prove this by an induction reminiscent of that
used in the proof of Theorem~\ref{thm:real-trans}.
In the course of their proof, they construct a polynomial $f$ with all roots
real having $d(k,n)$ (=degree of the Grassmannian) real points in 
$\pi_W^{-1}(f)$ (one for each chain in the
Bruhat order) and show that
$$
  \det d\pi_W\ =\ (-1)^{\omega(\ch)}\,,
$$
at the point in $\pi_W^{-1}(f)$ corresponding to the chain $\CH$.

White~\cite{Wh00} studied the statistic $|I(k,n)(-1)|$ and showed that it
equals zero if and only if $n$ is even, and that
$|I(2,2n)(-1)|=C_n$.
The main result of Eremenko and Gabrielov is the following, which gives
a lower bound for the number of real solutions in some cases of the Shapiro
conjecture. 

\begin{cor}[Eremenko and Gabrielov~{\cite[Corollary 2]{EG02}}]\label{cor:lbd}
 Let $L_1,L_2,\ldots,L_{k(n-k)}$ be codimension $(n-k)$-planes in
 ${\mathbb R}^n$ osculating the rational normal curve at $k(n-k)$ general real
 points. 
 Then number $\rho$ of real $k$-planes $K$ meeting each $L_i$
 non-trivially satisfies
$$
  |I(k,n)(-1)|\ \leq\ \rho\ \leq \ I(n,k)(1)\,.
$$
 In particular, when $n$ is odd this number $\rho$ is non-zero.
\end{cor}

As we remarked before, this lower bound holds if the osculating planes are no
longer required to be real, but that the set $\{t_1,t_2,\ldots,t_{k(n-k)}\}$ of
points of osculation are real (the roots of a real polynomial).

Consider now a linear projection $\pi_E$ different from the Wronski map.
If we move the center $E$, the degree of $\pi_E$ does not change as long as
$E$ does not meet the Grassmannian.
Thus if $E$ is in the same component of the space of codimension $k(n-k)+1$
planes not meeting the Grassmannian as is the center of the Wronski projection
$\pi_W$, then Corollary~\ref{cor:lbd} applies to $\pi_E^{-1}(x)$ for real $x$.

These new ideas of Eremenko and Gabrielov, particularly their notion of real
degree, greatly increase our understanding of the real Schubert calculus.
We emphasize that this is an important start in the search for lower bounds to
other enumerative problems. 
Not all enumerative problems involve a linear projection of a submanifold of
projective space.
In fact, of the other enumerative problems we have considered in the Schubert
calculus, only those of Theorem~\ref{thm:no-real} on the Lagrangian
Grassmannian involve linear projections, and they have a lower bound of 0.
There has been very little experimental work investigating 
lower bounds.
We describe some of it in the final section.

\subsection{Lower bounds in the Schubert calculus for flags?}
Conjecture~\ref{conj:flags} speculates that a zero-dimensional intersection of
Grassmannian Schubert varieties in a flag manifold $\Fla$ has only real
points, when the Schubert varieties are given by flags osculating the
rational normal curve at points and when the indices of the 
Grassmannian Schubet data are ordered.
If the indices cannot be ordered, then we conjectured (and have found
experimentally) that there is a selection of osculating flags with all points
of intersection real, as well as a selection with not all points real.

Table~\ref{table:12-flag} shows the results of computing 160,000 instances of the
 \begin{table}[htb]
  \begin{tabular}
   {|c||c|c|c|c|c|c|c|}\hline
   {Necklace} & \multicolumn{7}{c|}{Number of Real Solutions\rule{0pt}{13pt}}\\
   \cline{2-8}
        &0&2&4&6&8&10&12\rule{0pt}{13pt}\\\hline\hline
   22223333
     &0 & 0 & 0 & 0 & 0 & 0 & 20000\\\hline
   22322333
     &0 & 0 & 9 & 1677 & 3835 & 6247 & 8232\\ \hline
   22233233
     &0 & 0 & 67 & 3015 & 6683 & 4822 & 5413\\ \hline
   22332233
     &0 & 0 & 136 & 1533 & 7045 & 5261 & 6025\\ \hline
   22323323
     &0 & 0 & 303 & 2090 & 6014 & 7690 & 3903\\ \hline
   22323233
     &0 & 37 & 1944 & 4367 & 6160 & 4634 & 2858\\ \hline
   22232333
     &0 & 195 & 1476 & 1776 & 3628 & 4546 & 8379\\ \hline
   23232323
     &251 & 929 & 5740 & 3168 & 5420 & 2828 & 1662\\ \hline
  \end{tabular}
  \label{table:12-flag}
 \end{table}
%
%
intersection of 4 simple Schubert varieties $X_2\Fdot(t_i)$ and
 4 simple Schubert varieties $X_3\Fdot(t_i)$ in the manifold
 ${\mathbb F}\ell_{2<3}$  of partial flags 
 $E_2\subset E_3$ in 5-space.
Here the numbers $t_1<t_2<\cdots<t_8$ were 20,000  
random subsets of 8 numbers between 1 and 80, chosen using Maple's 
random number generator.
For each choice, we considered the 8 possibilities of orderings of the 
indices, and for each system we computed the number of real solutions 
out of 12.
Observe that the apparent lower bound on the number of real solutions depends
on ordering of the indices of the Schubert data. 
Other calculations we have done reinforce this observation.

We conclude by reminding the reader of the observation in Section 3 
that the number of real solutions to a problem in enumerative goemetry depends subtly on the configuration of the conditions.

\section*{Acknowledgements}
We thank the following people: Alexandre Eremenko, William Fulton, Andrei
Gabrielov, Viatcheslav Kharlamov, T.-Y. Li, Ricky Pollack, Felice Ronga,
Marie-Fran\c{c}oise Roy, Bernd Sturmfels, and  Thorsten Theobald. 

\def\cprime{$'$}

\end{document}